\documentclass[12pt]{article}
\usepackage{amsfonts}
\usepackage{psfrag}

\usepackage{graphicx}

\newtheorem{lemma}{\sc Lemma}[section]

\newtheorem{definition}[lemma]{\sc Definition}
\newtheorem{proposition}[lemma]{\sc Proposition}
\newtheorem{theorem}[lemma]{\sc Theorem}
\newtheorem{remark}[lemma]{\sc Remark}
\newtheorem{example}[lemma]{\sc Example}

\def\picture #1 by #2 (#3){
   \vbox to #2{
     \hrule width #1 height 0pt depth 0pt
     \vfill
     \special{picture #3} 
     }
   }

\def\scaledpicture #1 by #2 (#3 scaled #4){{
   \dimen0=#1 \dimen1=#2
   \divide\dimen0 by 1000 \multiply\dimen0 by #4
   \divide\dimen1 by 1000 \multiply\dimen1 by #4
   \picture \dimen0 by \dimen1 (#3 scaled #4)}
   }

\def\SS{\mathbb{S}}
\def\RR{\mathbb{R}}

\def\ii{\mathtt{i}}

\def\pp{\mathtt{p}}
\def\qq{\mathtt{q}}

\def\rr{\mathtt{r}}

\def\xx{\mathtt{x}}

\def\yy{\mathtt{y}}

\def\gotE{\mathfrak{E}}
\def\gotF{\mathfrak{F}}

\def\qed{\vrule height5pt depth0pt width5pt} 
\def\[{[\![}
\def\]{]\!]}
\def\ltriple{[\![\![}
\def\rtriple{]\!]\!]}

\def\osc#1#2{\big<#1\,\big\vert\,#2\big>}
\def\truc#1{\mathop{\triangleright}\limits_{\raise5pt\hbox{$_{#1}$}}}
\def\emptyword{{\circ\!\!\!/}}

\def\0{{\bf 0}}
\def\1{{\bf 1}}
\def\2{{\bf 2}}
\def\3{{\bf 3}}
\def\4{{\bf 4}}
\def\5{{\bf 5}}
\def\bfm{{\bf m}}

\begin{document}

\baselineskip=22pt


\centerline{\large\bf On the Gibbs properties of
Bernoulli convolutions}

\centerline{\large\bf related to
$\beta$-numeration in multinacci bases\footnote{This paper is dedicated to the
memory of Jean Marie Dumont.} }
\bigskip

\centerline{by {\sc Eric Olivier\footnote{Supported by a HK RGC
grant and a CUHK Postdoctoral Fellowship.}, Nikita
Sidorov\footnote{Supported by the EPSRC grant no GR/R61451/01.} \&
Alain Thomas} }
\bigskip
\centerline{---------------}
\bigskip
\centerline{\it (Release of the 11th of October 2004) }
\bigskip

{\leftskip=2cm
\rightskip2cm

{\small {\bf Abstract.} We consider infinitely convolved Bernoulli
measures (or simply Bernoulli convolutions) related to the
$\beta$-numeration. A matrix decomposition of these measures is
obtained in the case when $\beta$ is a PV number. We also
determine their Gibbs properties for $\beta$ being a multinacci
number, which makes the multifractal analysis of the corresponding
Bernoulli convolution possible. \par
\par}\par
}

\bigskip


\centerline{\large \bf 0. Introduction}
\medskip
The {\it Bernoulli convolutions} have been studied since the early
1930's \cite{JW, Erd, KS,Gar} and more recently, in the 1990's
onward, following the work of Alexander and Yorke \cite{AY}. They
have also been considered in view of their applications to fractal
geometry \cite{LP, LN, Lal, PS, KSS} and ergodic theory \cite{SV,
Sid} (see \cite{PSS} for further details and references). Our
approach is motivated by the fact that the Bernoulli convolution
associated with the golden ratio (usually called the {\em Erd\H os
measure}) proves to be weak Gibbs \cite{FO} and thus satisfies the
multifractal formalism. In the present paper we aim to generalize
this result.

\qquad We focus our attention on the Bernoulli convolutions
related to $\beta$-numeration for \hbox{$1<\beta<{\bf  d}$,} where
${\bf  d}\ge 2$ is an integer. The set ${\cal D}^{\bf N}$ of
sequences $\omega=(\omega_k)_{k=0}^\infty$ with {\it digits}
$\omega_k$ in the finite alphabet ${\cal D}:=\{\0,\dots,{\bf
d-1}\}$ is endowed with the product\hbox{ topology.} Given a ${\bf
d}$-dimensional probability vector $\pp:=(\pp_i)_{i=0}^{{\bf
d-1}}$, the {\em $\pp$-distributed $(\beta,{\bf  d})$-Bernoulli
convolution} is by definition the measure $\mu$ which corresponds
to the distribution of the random variable $X:{\cal D}^{\bf
N}\to{\bf R}$ such that
$X(\omega)=\sum\nolimits_{k=0}^\infty{\omega_k/ \beta^{k+1}}$,
where $\omega\mapsto\omega_k$ ($k=0,1,\dots$) is a sequence of
i.i.d. random variables, assuming the value
$i\in\{\0,\1,\dots,{\bf d-1}\}$ with the probability $\pp_i$. We
speak of the {\em uniform} Bernoulli convolution when
$\pp_i=1/{\bf  d}$ for each $i$.

\qquad The measure $\mu$ is self-similar and thus, satisfies the
so-called {\it law of pure types}. Recall that this means that
$\mu$ is either absolutely continuous or purely singular with
respect to the Lebesgue measure. Moreover, when each $\pp_i$ is
positive, the measure $\mu$ is fully supported by an interval.

\qquad Because of their nontrivial multifractal structure, we will
consider the Bernoulli convolutions known to be purely singular,
namely, those parameterized by Pisot-Vijayaraghavan (PV) number
$\beta$. In Section~2 we use an important arithmetic property of
PV numbers (Garsia's separation lemma) to obtain, in
Lemma~\ref{Ath4}, a decomposition of the Bernoulli measure $\mu$
which involves a matrix product (see \cite{Lal} for a similar
approach). Moreover, from the $\beta$-shift being of finite type
it follows (Section~2.3 and 2.4) that $\mu$ may be decomposed into
a finite number of measures having a specific structure; we call
them ${\cal M}$-measures, where ${\cal M}$ is a finite set of
square matrices with nonnegative entries.

\qquad In Section~3 we consider the case when $\mu$ is the
$(\beta,2)$-Bern\-oul\-li convolution, where $\beta$ is a {\em
multinacci number} of order $\bfm\ge2$ (see the definition below).
The main result of this section (Theorem~\ref{thA}) concerns the
``local" Gibbs properties of $\mu$.

\qquad Finally, in Section~4 we show that the ``local" Gibbs
properties allow one to apply the multifractal analysis of $\mu$
(see Theorem~\ref{thmultana}). Moreover, in the same section we
show how to use the multifractal analysis as a classification tool
and discuss the existence of ``global" Gibbs properties of $\mu$
in the cases of the uniform/nonuniform Erd\H os measure ($\bfm=2$
and $\beta=(1+\sqrt{5})/2$)---Theorems~\ref{erdosunif},
\ref{local2} and \ref{local3}.


\begin{section}{\large \bf Basic notions: definitions and
generalities}

{\bf 1.1. -- Net and adapted system of affine contractions --} For a given integer
${\bf s}\ge2$ we consider the finite alphabet
${\cal A}:=\{\0,\dots,{\bf s-1}\}$ and call any element $w=\xi_0\dots\xi_{n-1}
\in{\cal A}^n$
($n\ge1$) a {\em
word}. By convention, ${\cal A}^0:=\{\emptyword\}$, where
$\emptyword$ denotes the {\em empty word}. Let
${\cal A}^*$ denote the set of all words in the alphabet $\cal A$, {\em
i.e.}, ${\cal A}^*:=\bigcup_{n=0}^\infty{\cal A}^n$. The
concatenation of the two words $w,w'$ is, as usual, denoted by $ww'$.

\qquad The family $\gotF:=\{\gotF_n\}_{n=0}^\infty$ is called an
{\em ${\bf s}$-fold net} of the interval $[a,b]$ if $\gotF_n$ is a
partition of $[a,b[$ by ${\bf s}^n$ semi-open intervals (which we
denote by $\[w\]$ for $w\in{\cal A}^n$) with the extra property
that $\[ww'\]\subset\[w\]$, for any words $w,w\in{\cal A}^*$. Any
interval $\[w\]$ is by definition a {\it basic interval} of the
net $\gotF$. By Kolmogorov's Consistency Theorem, any positive
Borel measure $\eta$ on the real line whose support is a subset of
$[a,b]$, is characterized by its values taken on the intervals of
the ${\bf s}$-fold net.

\qquad Let ${\cal S}:=\{\SS_i\}_{i=0}^{{\bf s-1}}$  be a system of
(orientation preserving) affine contractions (s.a.c.) of the real
line; we say that ${\cal S}$ is {\it adapted} to the interval
$[a,b[$ when $\{\SS_i[a,b[\}_{i=0}^{{\bf s-1}}$ is a partition of
$[a,b[$. By convention, $S_{\emptyword}$ is the identity map on
$\bf R$ and, for any non-empty word $w=\xi_0\dots\xi_{n-1}\in{\cal
A}^*$, we put $S_w:=S_{\xi_0}\circ\cdots\circ S_{\xi_{n-1}}$. The
{\em ${\cal S}$-net} is by definition the ${\bf s}$-fold net
$\gotF:=\{\gotF_n\}_{n=1}^\infty$, where $\gotF_n$ is the
collection of the intervals $\[w\]:=\SS_w[a,b[$ for $w\in{\cal
A}^n$.
%

\qquad{\bf 1.2. -- ${\cal M}$-measures -- } Given ${\cal
M}:=\{M_i\}_{i=0}^{{\bf s-1}}$, a family of $r\times r$ matrices
($r\ge1$) with nonnegative entries, we denote
$M_*=M_0+\dots+M_{{\bf s-1}}$. Assume there exists a
column-vector~$R\not=0$ with nonnegative entries such that
$M_*R=R$; we say that $\nu$ is a {\em ${\cal M}$-measure} w.r.t.
the ${\bf s}$-fold net $\gotF$ if its support is a subset of
$[a,b]$ and if there exists a row-vector $L$ with nonnegative
entries such that $\nu\[w\]=LM_wR$, for any word
$w=\xi_0\dots\xi_{n-1}\in{\cal A}^*$ (by a similar convention as
above, $M_{\emptyword}$ is the unit $r\times r$ matrix and
$M_w:=M_{\xi_0}\cdots M_{\xi_{n-1}}$).
\medskip

\qquad {\bf 1.3. -- Gibbs, weak Gibbs and quasi-Bernoulli measures
--} Let $\eta$ be a finite positive Borel measure whose support is
$[a,b]$. We denote by $\sigma:{\cal A}^{\bf N}\to{\cal A}^{\bf N}$
the one-sided shift map defined for any
$\xi=(\xi_k)_{k=0}^\infty\in{\cal A}^{\bf N}$ by
$\sigma\xi=(\xi_{k+1})_{k=0}^\infty$. Recall that the {\em product
topology} on ${\cal A}^{\bf N}$ is  given  by the metric such that
the distance between $\xi$ and $\zeta$ is $2^{-k}$, where $k$ is
the length of the largest common prefix of $\xi$ and $\zeta$.
%

\begin{definition}{The measure $\eta$ is $\gotF$-weak Gibbs
in the sense of M. Yuri \cite{Yur}
if there exists a continuous map $\Phi$ from ${\cal
A}^{\bf N}$ to ${\bf R}$ (called a {\em potential}) and a
sequence of real numbers $K_n>1$, subexpo\-nential in the sense
that $\lim_n{1\over n}\log K_n=0$ such that, for any $\xi\in{\cal
A}^{\bf N}$ and any $n\ge1$,
\begin{equation}
{1\over K_n}\leq{\eta\[\xi_0\dots\xi_{n-1}\]\over
\exp\Big(\sum_{k=0}^{n-1} \Phi(\sigma^k\xi)\Big)}\leq
K_n.\label{A1}
\end{equation}
When a sequence $(K_n)$ can be taken constant, $\eta$ is {\em
$\gotF$-Gibbs} in the sense of Bowen \cite{Bow}.}
\end{definition}

%

For each $\xi\in{\cal A}^{\bf N}$ we put
$\phi_1(\xi)=\log\eta\[\xi_0\]$ and for $n\ge2$,
\begin{equation}\label{nstep}
\phi_n(\xi)=
\log\left({\eta\[\xi_0\cdots\xi_{n-1}\]\over\eta\[\xi_1\cdots\xi_{n-1}\]}\right).
\end{equation}
The continuous map $\phi_n:{\cal A}^{\bf N}\mapsto{\bf R}$
($n\ge1$) is called the {\em $n$-step potential} of $\eta$. Assume
that the sequence $\phi_n$ converges uniformly to a potential
$\Phi$; it is then straightforward that for~$n\ge1$,
\begin{equation}\label{A*}
{1\over K_n}\leq{\eta\[\xi_0\cdots\xi_{n-1}\]\over
\exp\Big(\sum_{k=0}^{n-1} \Phi(\sigma^k\xi)\Big)}\leq K_n
\quad\hbox{with}\quad K_n=\exp\left(\sum_{k=1}^n\Vert\Phi-\phi_n\Vert_\infty\right).
\end{equation}
By a well known lemma on the Ces\`aro sums, the sequence $(K_n)_n$
is subexponential, whence (\ref{A*}) means $\eta$ is $\gotF$-weak
Gibbs. Note that if $\sum_n\Vert\Phi-\phi_n\Vert_\infty<+\infty$,
then the $K_n$ are bounded and therefore, $\eta$ is a
$\gotF$-Gibbs measure.

\qquad We also consider quasi-Bernoulli measures.

\begin{definition}\label{quasibern}{\sl
The positive measure $\eta$ whose support is a subset of the
interval $[a,b]$, is said to be {\it $\gotF$-quasi-Bernoulli} if
there exists a constant $K>1$ such that for any words
$w,w'\in{\cal A}^*$,
$$
{1\over K}\eta\[w\]\eta\[w'\]\le\eta\[ww'\]\le K\eta\[w\]\eta\[w'\].
$$
The net $\gotF$ itself is said {\it quasi-Bernoulli} if the
Lebesgue measure (restricted to
$[a,b]$) is $\gotF$-quasi-Bernoulli.}
\end{definition}

Notice that a Gibbs measure is always quasi-Bernoulli.

%

\qquad {\bf 1.4. -- Multifractal analysis --} We need a number of
extra definitions. Recall that the {\it local dimension} at a
point $x$ which belongs to the support of the measure $\eta$, is
$\lim\limits_{r\to0}\log \eta\,\bigl(B_r(x)\bigr)/ \log
r$---provided the limit exists. (Here $B_r(x)$ stands for the
closed ball of radius $r$ centered at $x$.) Given an arbitrary
real $\alpha$, the {\it level set}~$E_\eta(\alpha)$ is defined as
the set of $x$ in the support of $\eta$ such that the local
dimension at $x$ exists and is equal to $\alpha$. The {\it
multifractal domain} $\hbox{\sc Dom}(\eta)$ is the set of
$\alpha\in{\bf R}$ for which $E_\eta(\alpha)$ is nonempty. The
{\it singularity spectrum} is the map which associates to any
$\alpha\in{\bf R}$ the Hausdorff dimension $\dim_HE_\eta(\alpha)$.
The {\it scale spectrum} (also called {\it $L^q$-spectrum}) is the
map $\tau_\eta$ from ${\bf R}$ to ${\bf R}\cup\big\{+\infty\big\}$
defined as follows:
$$
\tau_\eta(q):=\liminf_{r\to0}{\log
\inf\Big\{\sum_i\eta\big(J_i\big)^q\;;\;
\{J_i\}_i\Big\}\over \log r},
$$
where $\{J_i\}_i$ runs over the family of covers of the support of
$\eta$ by closed intervals  $J_i$ whose length is equal to $r$.
(We refer to the Book of Y. Pesin \cite{Pes} for analogue and
equivalent definitions of the scale spectrum.)

\qquad One usually says that $\eta$ satisfies the {\it
multifractal formalism} if the singularity spectrum and the scale
spectrum of $\eta$ form a {\it Legendre transform pair}. The
multifractal formalism is trivially not universal, but it has been
established for wide classes of measures \cite{Ran,BMP,Hof,PW,Ols}, when
some conditions of geometric homogeneity are satisfied
(self-similarity, conformality or Gibbs properties, etc.). For our
purpose, we refer to the multifractal formalism of the
quasi-Bernoulli and weak Gibbs measures stated in the two
following theorems:

\begin{theorem}\label{herteau}
{\sl\cite{BMP, Heu} Let $\gotE$ be an ${\bf s}$-fold
quasi-Bernoulli net of a compact interval $[a,b]$ and $\eta$ be a
positive measure fully supported by $[a,b]$. If $\eta$ is a
$\gotE$-quasi-Bernoulli measure, then
\begin{enumerate}
\item The scale spectrum $\tau_\eta$ of $\eta$ is concave and
differentiable on the whole real line and moreover,
$$
-\infty<\underline\alpha:=\lim_{q\to+\infty}{\tau_\eta(q)\over
q}\le \lim_{q\to-\infty}{\tau_\eta(q)\over
q}=:\overline\alpha<+\infty.
$$
\item The multifractal domain of $\eta$ is $\hbox{\sc
Dom}(\eta)=[\underline{\alpha};\overline{\alpha}]$ and for any
$\underline{\alpha}\le \alpha\le \overline{\alpha}$,
$$
\dim_HE_\eta(\alpha)=\inf_{q\in{\bf R}}\{\alpha q-\tau_\eta(q)\},
$$
meaning that $\eta$ satisfies the multifractal formalism.
\end{enumerate}}
\end{theorem}

\begin{theorem}\label{Fengolivier}
{\sl\cite[Theorem A$'$]{FO} Let $\gotE$ be an ${\bf s}$-fold net
of a compact interval $[a,b]$ with respect to which the Lebesgue
measure is Gibbs. If the positive measure $\eta$, supposed to be
fully supported by $[a,b]$, is weak Gibbs w.r.t. $\gotE$, then
\begin{enumerate}
\item The scale spectrum $\tau_\eta$ of $\eta$ is concave on the
whole real line and moreover,
$$
-\infty<\underline\alpha:=\lim_{q\to+\infty}{\tau_\eta(q)\over
q}\le \lim_{q\to-\infty}{\tau_\eta(q)\over
q}=:\overline\alpha<+\infty.
$$
\item The multifractal domain of $\eta$ is $\hbox{\sc
Dom}(\eta)=[\underline{\alpha};\overline{\alpha}]$ and for any
$\underline{\alpha}\le \alpha\le \overline{\alpha}$,
$$
\dim_HE_\eta(\alpha)=\inf_{q\in{\bf R}}\{\alpha q-\tau_\eta(q)\},
$$
meaning that $\eta$ satisfies the multifractal formalism.
\end{enumerate}
}
\end{theorem}

\begin{remark}\label{classifyingtool}{The analysis
of the Gibbs properties of a given measure is a simple way to
study its multifractal structure. Conversely, the multifractal
properties of the  measure may be used as a classification tool
w.r.t. its Gibbs properties.

\qquad For instance, a Gibbs measure is both weak Gibbs and
quasi-Bernoulli, but a weak Gibbs measure need not be
quasi-Bernoulli, for there exist weak Gibbs measures, whose scale
spectrum is not differentiable (see \cite{Fen,FO}).

\qquad Another interesting application uses topological properties
of the multifractal domain for the measure in question. Actually,
a measure whose multifractal domain is disconnected (or
noncompact) is neither weak Gibbs nor quasi-Bernoulli w.r.t. any
given reasonable net.

\qquad Below we discuss these types of classification for the
cases of the uniform/nonuniform Erd\H os measure---see
Theorem~\ref{local2} and Theorem~\ref{local3}.}
\end{remark}

\end{section}

\begin{section}{Decomposition of Bernoulli convolutions}

{\bf 2.1. -- Generalities --} Let ${\bf b},{\bf d}$ be two
integers and $\beta$ a real number such that
$$
1\le{\bf b}-1<\beta\le {\bf b}\le{\bf d}.
$$
Throughout this paper, we consider the following two alphabets:
$$
{\cal B}:=~\big\{\0,\dots,{\bf b}-1\big\}\quad\hbox{and}\quad{\cal
D}:=\big\{\0,\dots,{\bf d}-1\big\}.
$$
Recall that the $\pp$-distributed $(\beta,{\bf d})$-Bernoulli convolution is
the measure $\mu$ such that, for any Borel set $B\subset{\bf R}$,
$$
\mu(B)=I\!\!P\{\omega\; : \;X(\omega)\in B\},
$$
where $X(\omega)=\sum\nolimits_{k=0}^\infty{\omega_k/ \beta^{k+1}}$
for any $\omega\in{\cal D}^{\bf N}$, and $I\!\!P$ is
the product measure on ${\cal D}^{\bf N}$ with
parameter $\pp=\left(\pp_0,\dots,\pp_{{\bf d}-1}\right)$.

\qquad The $\beta$-numeration is usually related to the s.a.c.
$\{\RR_i\}_{i=0}^{{\bf b}-1}$ defined as follows:
$$
\RR_i(x)=(x+i)/\beta.
$$
From now on, $i$ will always stand for an arbitrary fixed element
in ${\cal B}$; denoting by $\sigma:{\cal D}^{\bf N}\to{\cal
D}^{\bf N}$ the shift map on the product space ${\cal D}^{\bf N}$,
for any Borel set $B\subset{\bf R}$, and $x\in{\bf R}$, one has
$X(\omega)\in B+x$ if and only if
$X(\sigma\omega)\in\RR_i^{-1}(B)+\beta x+(i-\omega_0)$. Since
$I\!\!P$ is a product measure,
\begin{equation}\label{A0}
\mu(B+x)=\sum_{j={\bf 0}}^{{\bf d}-1}\pp_j\cdot
\mu\Big(\RR_i^{-1}(B)+\beta
x+(i-j)\Big).
\end{equation}
Suppose that $\RR_i^{-1}(B)\subset[0,1]$; the support of $\mu$ being a subset of $[0,\alpha_\mu]$
with
$\alpha_\mu:=({\bf d}-1)/(\beta-1)$, one has $\mu(\RR_i^{-1}(B)+y)=0$ whenever
$y\notin]-1,\alpha_\mu[$. Given any $x\in]-1,\alpha_\mu[$, we
write,
$$
x\mathop{\triangleright}\limits^i y\iff -1<y<
\alpha_\mu\quad\hbox{and}\quad\beta x+i-y\in{\cal D}
$$
and
$$
x\triangleright y\iff  \exists i\in{\cal B},\;x\mathop{\triangleright}\limits^i y.
$$
Then, in order to simplify the relation (\ref{A0}), we set the
following definition:
%

\begin{definition}\label{DefI}
{\sl The real number $r$ belongs to the set
${\cal I}_{(\beta,{\bf d})}$ if there exist
a finite sequence of real numbers $0=r_0,\dots,r_{n}=r$ such that
$r_0\triangleright r_1\triangleright \cdots
\triangleright r_{n}$.}
\end{definition}

%
The set ${\cal I}_{(\beta,{\bf d})}$ is always nonempty, because
it contains at least $0$; moreover, it  is clearly finite or
countable and we denote by $0=\ii_0,\ii_1,\dots$ the sequence of
its elements. Actually, in order to determine  ${\cal
I}_{(\beta,{\bf d})}$ explicitly, we present an equivalent
definition using induction (see also \cite{Lau} and \cite[Section
2]{BH}).

\qquad Firstly, since $0=\ii_0\in {\cal I}_{(\beta,{\bf d})}$, we define
${\cal I}_0:=\{{\ii_0}\}$; then, assuming that ${\cal
I}_n=\{\ii_0,\dots,\ii_{k_n}\}$, we put
$$
{\cal
I}_{n+1}:=\bigcup_{k=0}^{k_n}\{y=\beta\ii_k+(i-j)\;;\;(i,j)\in{\cal
B}\times{\cal D}\;\hbox{and}\;-1<y<\alpha_\mu\},
$$
and finally,
\begin{equation}\label{defibetad}
{\cal I}_{(\beta,{\bf d})}:=\bigcup_n{\cal I}_n.
\end{equation}
Notice that ${\cal I}_{(\beta,{\bf d})}$ is finite  whenever
${\cal I}_{n-1}={\cal I}_n$ for some $n\ge1$; in that case, ${\cal
I}_{(\beta,{\bf d})}={\cal I}_n$.

\qquad It follows from (\ref{A0}) that,
for $B\subset[0,1]$ with $\RR_i^{-1}(B)\subset[0,1]$,
\begin{equation}
\mu(B+\ii_h)=
\sum\nolimits_kM_i(h,k)\,\mu(\RR_i^{-1}(B)+\ii_k)
\end{equation}
where $M_i$ is the nonnegative matrix defined as follows: for the
row index $h$ and the column index $k$, with $\ii_h$ and $\ii_k$
in ${\cal I}_{(\beta,{\bf d})}$, by definition,
\begin{equation}\label{A9}
M_i(h,k)= \cases{ \pp_j, & if $j=i+\beta \ii_h-\ii_k\in{\cal
D}$\cr\cr 0,&otherwise.\cr}
\end{equation}
%

\begin{lemma}
\label{Ath4} {\sl Let $i\in{\cal B}$ and assume ${\cal
I}_{(\beta,{\bf d})}:= \big\{\ii_0,\dots,\ii_{{\bf r}-1}\big\}$;
then, for any Borel set $B\subset[0,1]$ with $\RR_i^{-1}(B)\subset[0,1]$, one has:
$$
\pmatrix{
\mu(B+\ii_0)\cr
\vdots\cr
\mu(B+\ii_{{\bf r}-1})}
=M_{i
}
\pmatrix{
\mu\Big(\RR_i^{-1}(B)+\ii_0\Big)\cr
\vdots\cr
\mu\Big(\RR_i^{-1}(B)+\ii_{{\bf r}-1}\Big)\cr}.
$$}
\end{lemma}

%
\begin{figure}[p]
  \begin{center}
       \includegraphics{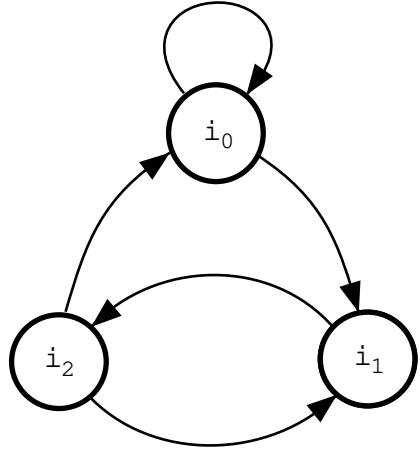}
       \caption{\leftskip=1.5cm\rightskip=1.5cm
{ \small\it In the case of the PV number $\beta$ such that $\beta^2=5\beta+3$, the set ${\cal
I}_{(\beta,6)}$ has three elements, say $\ii_0=0$, $\ii_1=1$ and $\ii_2=\beta-5$; here, we give a
representation of the graph of the relation $\cdot\triangleright\cdot$ on ${\cal
I}_{(\beta,6)}$.}}
  \end{center}
\end{figure}
%
\begin{figure}[p]
  \begin{center}
       \psfrag{X}[][]{Automaton of $\cdot \mathop{\triangleright}\limits^i \cdot$ with $i={\bf 0}$}
       \psfrag{Y}[][]{Automaton of $\cdot \mathop{\triangleright}\limits^i \cdot$ with $i={\bf 1}$}
       \psfrag{Z}[][]{Automaton of $\cdot \mathop{\triangleright}\limits^i \cdot$ with $i={\bf 2}$}
       \psfrag{T}[][]{Automaton of $\cdot \mathop{\triangleright}\limits^i \cdot$ with $i={\bf 3}$}
       \psfrag{U}[][]{Automaton of $\cdot \mathop{\triangleright}\limits^i \cdot$ with $i={\bf 4}$}
       \psfrag{V}[][]{Automaton of $\cdot \mathop{\triangleright}\limits^i \cdot$ with $i={\bf 5}$}
       \includegraphics{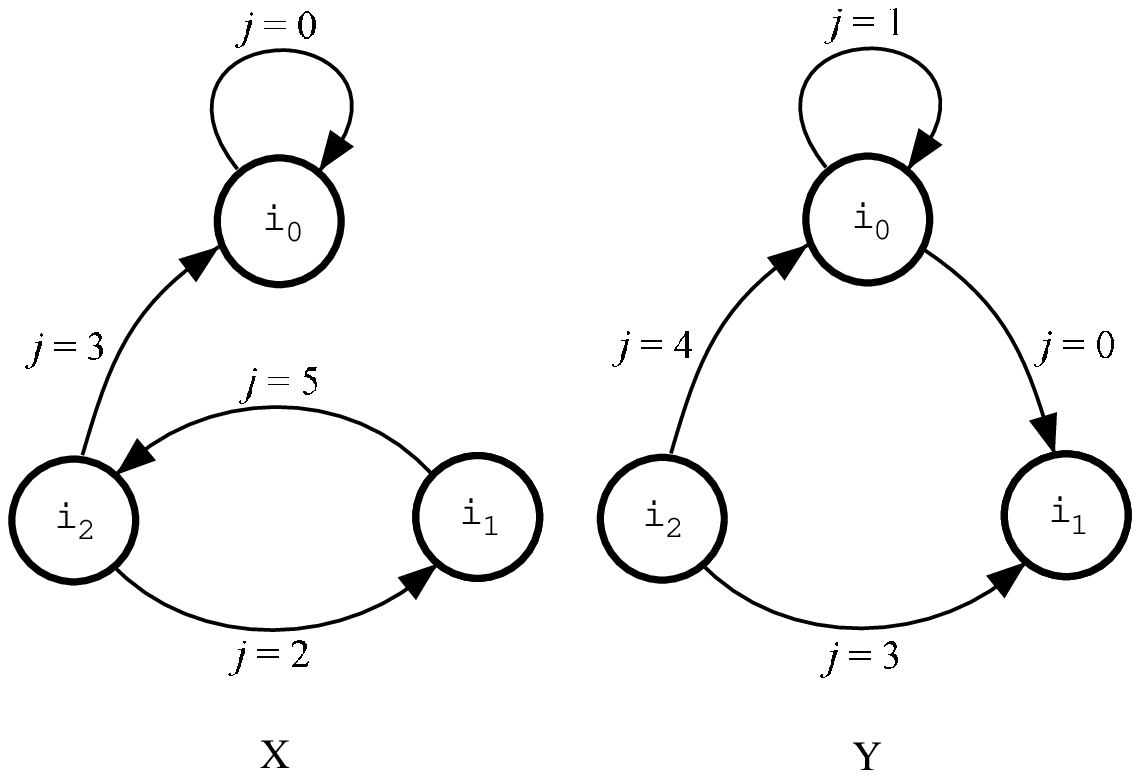}
       \caption{\leftskip=1.5cm\rightskip=1.5cm
{ \small\it We represent here, the automaton of the relation $\cdot
\mathop{\triangleright}\limits^i
\cdot$ for $i={\bf 0}$ and ${\bf 1}$ (the cases $i={\bf 2},\dots,{\bf 5}$ may be recovered by mean of the
matrices in (\ref{mat})). When
$\ii_h\mathop{\triangleright}\limits^i\ii_k$,  the label $j$ of the arrow from state $\ii_h$ to state $\ii_k$,
means that
$j=\beta\ii_h-\ii_k+i$; then, $j\in{\cal D}$ and according to (\ref{A9}) one has
$M_i(h,k)=\pp_j$.}}
  \end{center}
\end{figure}

\begin{example}
{\sl Let $\beta$ be the algebraic integer such that
$\beta^2=5\beta+3$ with $5<\beta<6$ and ${\bf d}=6$, so that
${\cal B}={\cal D}=\{\0,\dots,\5\}$. We apply the argument
described above to determine the set ${\cal I}_{(\beta,6)}$ and
the matrices $M_{\bf 0},\dots,M_{\bf 5}$. Following the induction
process leading to the definition in (\ref{defibetad}), one has
${\cal I}_0=\{\ii_0=0\}$ and using the fact that
$1<\alpha_\mu=5/(\beta-1)<2$, one obtains successively:
\begin{eqnarray*}
&&{\cal I}_1={\cal I}_0\cup\{y=i-j\;;\;(i,j)\in{\cal B}\times{\cal
D}\;\hbox{and}\;-1<y<\alpha_\mu\}=\{\ii_0=0,\ii_1=1\}\\
&&\qquad\hbox{(\it with $0\triangleright0$ and
$0\triangleright1$)}\\
&&{\cal I}_2={\cal I}_1\cup\{y=\beta+(i-j)\;;\;(i,j)\in{\cal B}\times{\cal
D}\;\hbox{and}\;-1<y<\alpha_\mu\}\\
&&\;\quad =\{\ii_0=0,\ii_1=1,\ii_2=\beta-5\}\\
&&\qquad \hbox{(\it with,
$1\triangleright\beta-5$)}\\
&&{\cal I}_3={\cal I}_2\cup\{y=\beta(\beta-5)+(i-j)\;;\;(i,j)\in{\cal B}\times{\cal
D}\;\hbox{and}\;-1<y<\alpha_\mu\}={\cal I}_2={\cal I}_{(\beta,6)}\\
&&\qquad\hbox{(\it with
$\beta-5\triangleright0$
and $\beta-5\triangleright1$)}
\end{eqnarray*}
Finally, $ {\cal
I}_{(\beta,6)}=\{\ii_0=0,\ii_1=1,\ii_2=\beta-5\}$. {\it Figure~1}
shows the graph of the relation~$\cdot\triangleright\cdot$ on
${\cal I}_{(\beta,6)}$. Moreover, the algorithm we use to
determine the set ${\cal I}_{(\beta,6)}$ provides us with extra
information sufficient to obtain the matrices $M_{\bf
0},\dots,M_{\bf 5}$ as defined in (\ref{A9})---see {\it Figure~2}.
We thus have (for $i=\0,\dots,\5$):
\begin{equation}\label{mat}
M_i=\pmatrix{\pp_i&\pp_{i-1}&0\cr0&0&\pp_{i+5}\cr\pp_{i+3}&\pp_{i+2}&0},
\end{equation}
where, by convention, $\pp_i=0$,  for any $i\le-1$ or $i\ge6$ (see
\cite{OST4} for the  general case when $\beta$ is a quadratic
number). }
\end{example}

\bigskip
\bigskip

\begin{example}
{\sl Assume now $\beta$ to be the PV number satisfying
$\beta^3=3\beta^2-1$. The set ${\cal I}_{(\beta,3)}$ has then
eight elements which we present in the form of the list:
\begin{equation}\label{list}
\matrix{
\ii_0=0,\hfill&\ii_1=1,\hfill&\ii_2=\beta-2,\hfill&\ii_3=\beta^2-2\beta-2,\hfill\cr\cr
\ii_4=\beta^2-2\beta-3,\hfill&\ii_5=\beta^2-3\beta,\hfill&\ii_6=\beta^2-3\beta+1,\hfill&\ii_7=\beta-3.\hfill\cr}
\end{equation}

The graph of the relation
$\cdot\triangleright\cdot$ on ${\cal I}_{(\beta,3)}$, is represented in {\it Figure 3}. 
%
\begin{figure}[!h]
  \begin{center}
       \includegraphics{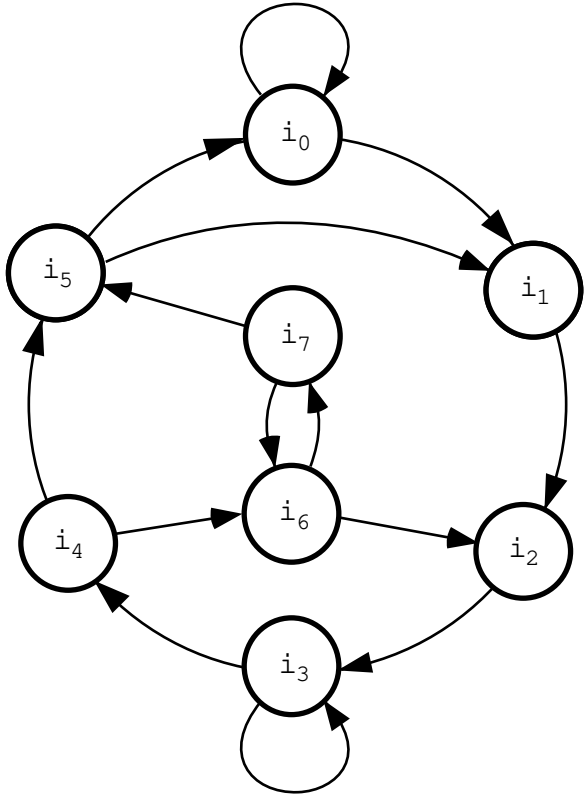}
       \caption{\leftskip=1.5cm\rightskip=1.5cm
{ \small\it In the case of the PV number $\beta$ such that
$\beta^3=3\beta^2-1$, the set ${\cal I}_{(\beta,3)}$ has eight
elements listed in (\ref{list}); here, we give a representation of
the graph of the relation $\cdot\triangleright\cdot$ on ${\cal
I}_{(\beta,3)}$.}}
  \end{center}
\end{figure}
%

\vbox{The 
matrices associated by (\ref{A9}) to the $(\pp\;\qq\;\rr)$-distributed $(\beta,3)$-Bernoulli
convolution are:
\begin{eqnarray*}
&&M_{\bf 0}=
\pmatrix{\pp&0&0&0&0&0&0&0\cr
0&0&\rr&0&0&0&0&0\cr
0&0&0&\rr&0&0&0&0\cr
0&0&0&\qq&\rr&0&0&0\cr
0&0&0&0&0&0&0&0\cr
0&0&0&0&0&0&0&0\cr
0&0&\qq&0&0&0&0&\rr\cr
0&0&0&0&0&\pp&0&0\cr
}
\qquad 
M_{\bf 1}=
\pmatrix{\qq&\pp&0&0&0&0&0&0\cr
0&0&0&0&0&0&0&0\cr
0&0&0&0&0&0&0&0\cr
0&0&0&\rr&0&0&0&0\cr
0&0&0&0&0&\pp&0&0\cr
\pp&0&0&0&0&0&0&0\cr
0&0&\rr&0&0&0&0&0\cr
0&0&0&0&0&\qq&\pp&0\cr}\\
&&M_{\bf 2}=
\pmatrix{\rr&\qq&0&0&0&0&0&0\cr
0&0&0&0&0&0&0&0\cr
0&0&0&0&0&0&0&0\cr
0&0&0&0&0&0&0&0\cr
0&0&0&0&0&\qq&\pp&0\cr
\qq&\pp&0&0&0&0&0&0\cr
0&0&0&0&0&0&0&0\cr
0&0&0&0&0&\rr&\qq&0\cr}.
\end{eqnarray*}
}
}
\end{example}

\qquad {\bf 2.2.-- Known cases when ${\cal I}_{(\beta,{\bf d})}$ is finite --} Our analysis of a $(\beta,{\bf
d})$-Bernoulli convolution relies on the finiteness of ${\cal I}_{(\beta,{\bf d})}$; for instance, this is trivially the
case  when
$\beta$ is an integer:
%

\begin{proposition}
{\sl If $\beta$ is an integer one has
$
{\cal I}_{(\beta,{\bf d})}=\{\0,\dots,{\bf a}-1\}
$,
where ${\bf
a}$ is the integer such that ${\bf a}-1<{({\bf d}-1)/(\beta-1)}\le{\bf a}$.}
\end{proposition}

%

This proposition is the starting point of our analysis of the
Bernoulli convolution in an integral basis developed in
\cite{OST5}; however, the methods used in that paper are different
from the ones we use here.

\qquad The second case for which we know that ${\cal
I}_{(\beta,{\bf  d})}$ is finite arises when $\beta$ is a
noninteger PV number, i.e., an algebraic integer $\beta=\beta_s>1$
whose Galois conjugates $\beta_1,\dots,\beta_{s-1}$ are strictly
less than 1 in modulus. Given any polynomial $A(X)\in{\bf Z}[X]$
such that $A (\beta)\not=0$, it is necessary that
$A(\beta_k)\not=0$ for $k=1,\dots,s$ and thus the integer $\vert
A(\beta_1)\cdots A(\beta_s)\vert$ is greater or equal to $1$.
Since $\vert\beta_k\vert<~1$ for $k=1,\dots,s-1$, one deduces that
$\vert A(\beta_k)\vert\le M/(1-\vert\beta_k\vert)$, where $M$ is
the maximum of the absolute values of the coefficients of $A(X)$:
in other words,
\begin{equation}\label{PV}
\forall A(X)\in{\bf Z}[X],\quad A(\beta)\not=0\;\Rightarrow\;\vert
A(\beta)\vert\ge
{1\over M^{s-1}}\prod_{k=1}^{s-1}(1-\vert\beta_k\vert)>0.
\end{equation}
The fact that (\ref{PV}) is satisfied by the PV numbers has been
discovered by Garsia \cite{Gar} and is usually called {\em
Garsia's separation lemma} (see also \cite{Lal}).

\qquad Return to the question of the cardinality of ${\cal
I}_{(\beta,{\bf  d})}$ and note that for any element $i$ in this set
there exists a finite sequence $\varepsilon_0,\dots,\varepsilon_m$ of
integers lying between $-{\bf  d}$ and ${\bf b}$ such that
$i=\varepsilon_0+\varepsilon_1\beta+\dots+\varepsilon_m\beta^m$. Since
$\beta$ is a PV number, it follows from (\ref{PV}) that the
distance between two different elements in ${\cal I}_{(\beta,{\bf  d})}$
is bounded from below by
$(2{\bf d})^{1-s}\prod_{k=1}^{s-1}(1-\vert\beta_k\vert)$, which
yields the following proposition:
%

\begin{proposition}
\label{Pisot} {\sl The set ${\cal I}_{(\beta,{\bf  d})}$ is finite
when $\beta$ is a PV number.}
\end{proposition}

\begin{proposition}If ${\cal I}_{(\beta,{\bf  d})}$ is finite,
then all the conjugates of $\beta$ are less than $\beta$ in
modulus and also less in modulus than $(1+\sqrt5)/2$.
\end{proposition}

The proof immediately follows from the fact that ${\cal
I}_{(\beta,{\bf  d})}$ contains the orbit of 1 under the
$\beta$-shift; then we may use the results of Solomyak \cite{Sol}.
We leave details to the reader.

%
\begin{remark}It should be interesting find a non-PV number $\beta$
for which ${\cal I}_{(\beta,{\bf d})}$ is finite, even when
(\ref{PV}) does not hold. We believe the only possible
counterexample may be a Salem number (the one for which all its
conjugates are less than or equal to 1 in modulus and some of them
are equal to 1 in modulus).
\end{remark}

\qquad {\bf 2.3.-- Heuristics --} In this section we give a
heuristic description of the framework that we use in Section~3
for the complete analysis of the Bernoulli convolution in the
multinacci bases.

\qquad Recall that a $(\beta,{\bf d})$-Bernoulli convolution $\mu$
is supported by the interval $[0,\alpha_\mu]$ with
$\alpha_\mu=({\bf d}-1)/(\beta-1)$. We assume that ${\cal
I}_{(\beta,{\bf d})}=\{\ii_0,\dots,\ii_{\rr-1}\}$ and for any
$j=0,\dots,{\bf r}-1$, we define the measure $\mu_j$ by putting ,
for any Borel set $B$ of the real line,
\begin{equation}\label{defmuj}
\mu_j(B)=\mu\Big(B\cap[0,1]+\ii_{j}\Big).
\end{equation}
If both $B$ and $\RR_i^{-1}(B)$ ($i\in{\cal B}$) are subsets of
$[0,1]$, it follows from Lemma~\ref{Ath4} that
\begin{equation}\label{correct0}
\pmatrix{
\mu_0(B)\cr
\vdots\cr
\mu_{{\bf r}-1}(B)\cr}=M_i
\pmatrix{
\mu_0\Big(\RR_i^{-1}(B)\Big)\cr
\vdots\cr
\mu_{{\bf r}-1}\Big(\RR_i^{-1}(B)\Big)\cr}.
\end{equation}
\begin{figure}[p]
  \begin{center}
       \psfrag{A}[][]{$\RR_1(x)\!=\!{x/\beta}\!+\!{1/ \beta}$}
       \psfrag{B}[][]{${1/ \beta}$}
       \psfrag{C}[][]{${1/\beta^2}$}
       \psfrag{D}[][]{$\RR_0(x)\!=\!{x/ \beta}$}
       \psfrag{E}[][]{$x$}
       \psfrag{X}[][]{$\widehat{\RR}_2(x)\!=\!\RR_{10}(x)$}
       \psfrag{Z}[][]{$\widehat{\RR}_1(x)\!=\!\RR_{010}(x)$}
       \psfrag{Y}[][]{$\widehat{\RR}_0(x)\!=\!\RR_{00}(x)$}
       \includegraphics{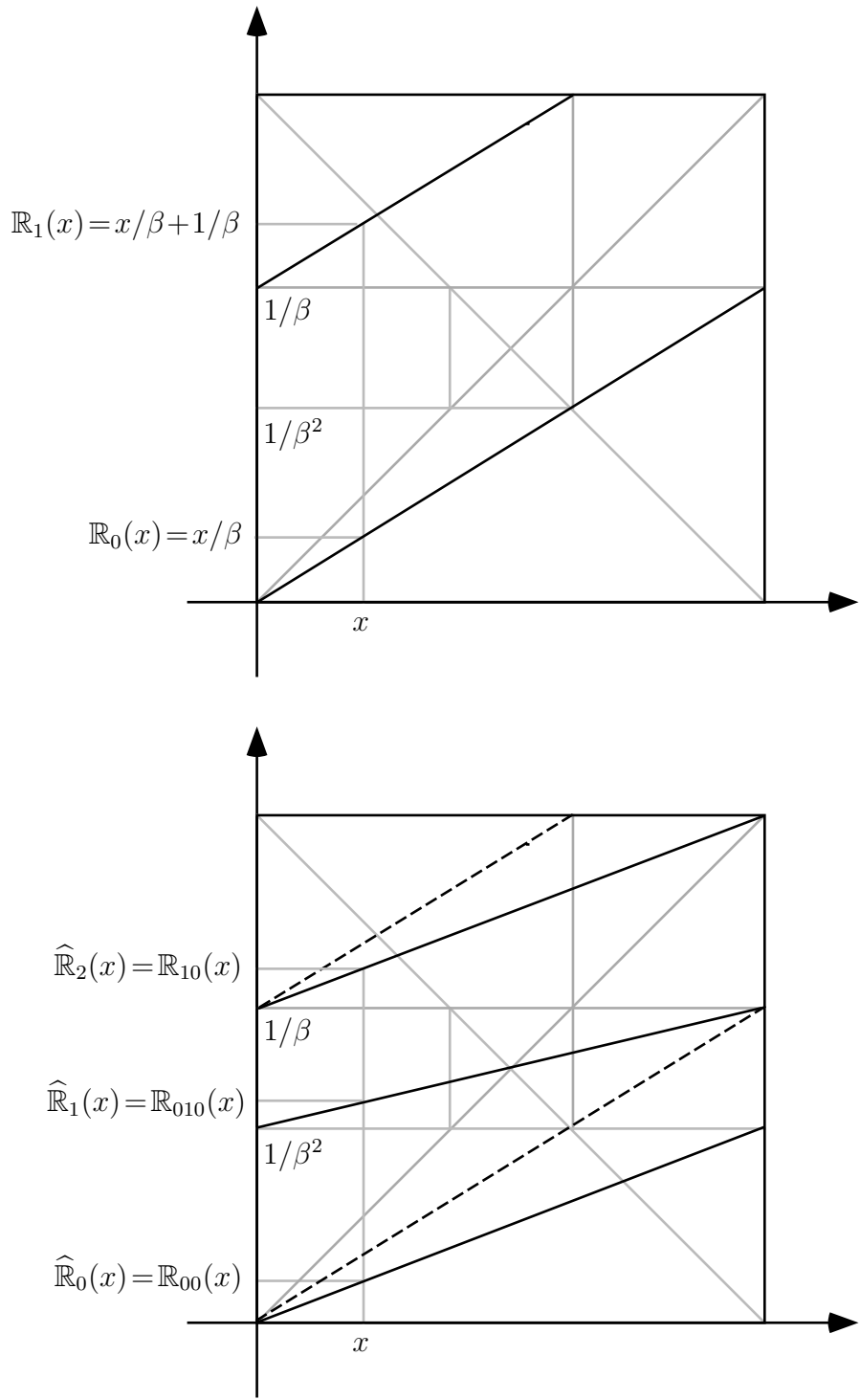}
       \caption{\leftskip=1.5cm\rightskip=1.5cm
{ \small\it The first plot contains the two affine contractions
$\RR_0$ and $\RR_1$ associated with the $\beta$-shift in the case
when $\beta=(1+\sqrt{5})/2$. In the second plot we give an example
of a s.a.c. ${\cal R}=\{\widehat{\RR}_k\}_{k=0}^3$ which is
adapted to the unit interval.}}
  \end{center}
\end{figure}
\qquad Our second assumption is the existence of ${\bf s}$ words
$w_{\bf 0},\dots,w_{{\bf s-1}}$ in ${\cal B}^*$ which satisfy the
two conditions: firstly, $\RR_w[0,1]\subset[0,1]$, for any suffix
$w$ of any of the words $w_j$, for $j={\bf 0},\dots,{\bf s-1}$;
secondly, the s.a.c. ${\cal
R}:=\{\widehat{\RR}_j=\RR_{w_j}\}_{j={\bf 0}}^{{\bf s}-1}$ is
adapted to the interval $[0,1]$.

\qquad The s.a.c ${\cal R}$ is associated to a ${\bf s}$-fold net of
$[0,1[$, say $\gotF$, whose basic intervals are coded by the words
in $\{\0,\dots,{\bf s-1}\}^*$. Recall that the basic interval of
$\gotF$ generated by a word
$\xi_0\cdots\xi_{n-1}\in\{\0,\dots,{\bf s-1}\}^n$ is by definition
$\[\xi_0\cdots\xi_{n-1}\]:=\widehat{\RR}_{\xi_0\cdots\xi_{n-1}}[0,1[$;
then one has successively:
\begin{eqnarray*}
&&\widehat{\RR}_{\xi_0}^{-1}\[\xi_0\cdots\xi_{n-1}\]=\[\xi_1\cdots\xi_{n-1}\]\subset[0,1],\\
&&\widehat{\RR}_{\xi_1}^{-1}\[\xi_1\cdots\xi_{n-1}\]=\[\xi_2\cdots\xi_{n-1}\]\subset[0,1],\\
&&\quad\vdots\\
&&\widehat{\RR}_{\xi_{n-1}}^{-1}\[\xi_{n-1}\]=[0,1[\subset[0,1].
\end{eqnarray*}
Therefore, denoting by $\widehat{M}_j:=M_{w_j}$, a recursive
application of (\ref{correct0}) yields:
\begin{equation}\label{correct1}
\pmatrix{ \mu_{0}\[\xi_0\cdots\xi_{n-1}\]\cr \vdots\cr \mu_{{{\bf
r}-1}}\[\xi_0\cdots\xi_{n-1}\]\cr}=\widehat{M}_{\xi_0\cdots\xi_{n-1}}
\pmatrix{ \mu_{0}[0,1[\cr \vdots\cr \mu_{{{\bf r}-1}}[0,1[\cr}.
\end{equation}
For ${\cal M}:=\{\widehat{M}_j\}_{j={\bf 0}}^{{\bf s}-1}$, it is
clear from (\ref{correct1}) that $\mu_0,\dots,\mu_{{\bf r}-1}$ are
${\cal M}$-measures w.r.t. the ${\bf s}$-fold net $\gotF$ .

\qquad Finally, we would like to make a simple remark which may
simplify the analysis of $\mu$. Let $U_0,\dots, U_{\rr-1}$ be the
elements of the canonical basis of the $1\times\rr$-matrix vector
space. We make the third assumption that any of the matrices
$\widehat{M}_j$ ($j={\bf 0},\dots,{\bf s-1}$)  leaves the vector
space generated by $U_{{k_1}},\dots,U_{{k_t}}$, invariant; we
denote by $\widehat{M}_j'$ ($j={\bf 0},\dots,{\bf s-1}$) the
corresponding $t$-dimensional submatrices. Then (\ref{correct1})
can be reduced to
\begin{equation}\label{correct2}
\pmatrix{ \mu_{k_1}\[\xi_0\cdots\xi_{n-1}\]\cr \vdots\cr
\mu_{k_t}\[\xi_0\cdots\xi_{n-1}\]\cr}=\widehat{M}_{\xi_0\cdots\xi_{n-1}}'
\pmatrix{ \mu_{k_1}[0,1[\cr \vdots\cr \mu_{k_t}[0,1[\cr}.
\end{equation}

\qquad Let us introduce an auxiliary measure which we denote by
$\mu_*$. It is associated to the Bernoulli convolution $\mu$ by
putting, for any Borel subset $B$ of the real line:
$$
\mu_*(B)={ \sum_{j=1}^t\mu_{k_j}(B)\over
\sum_{j=1}^t\mu_{k_j}[0,1[},
$$
so that, for any word $w\in\{\0,\dots,{\bf s}-1\}^*$,
\begin{equation}\label{correct3}
\mu_*\[w\]:=L\widehat{M}_w'R,
\end{equation}
where
$$
L:=\pmatrix{1\dots1\cr}\quad\hbox{and}\quad R:={1\over
\sum_j\mu_{k_j}[0,1[}\pmatrix{ \mu_{k_1}[0,1[\cr \vdots\cr
\mu_{k_t}[0,1[\cr}.
$$
\qquad In general, the Bernoulli convolution $\mu$ itself needs
not be an ${\cal M}$-measure, which is why one of the main
advantages of introducing successively the measures
$\mu_{k_1},\dots,\mu_{k_t}$ and $\mu_*$ is the fact that they all
are ${\cal M}$-measures. The measure~$\mu_*$ proves to be a better
candidate for studying its Gibbs properties; roughly speaking,
this is due to the fact that the left row vector in
(\ref{correct3}) has strictly positive entries. Thus, our study of
the Gibbs properties of $\mu_*$ (and, consequently, of $\mu$) will
be reduced to the analysis of the convergence of the $n$-step
potential  $\phi_n$ ($n=1,2\dots$)---which, in the present case,
is defined for any $\xi\in\{\0,\dots,{\bf s-1}\}^{\bf N}$~by
$$
\phi_n(\xi):=
\log\left({L\widehat{M}_{\xi_0\cdots\xi_{n-1}}'R\over
L\widehat{M}_{\xi_1\cdots\xi_{n-1}}'R}\right).
$$

\qquad As we will see, when $\mu$ is a Bernoulli convolution
associated with the multinacci numbers, the measure $\mu_*$
displays a clear Gibbs structure which will be analyzed in detail
below (Theorem~\ref{thA}) and which ensures that the multifractal
formalism holds. Nevertheless, it is worth noting that the
Bernoulli convolution $\mu$ itself may not be Gibbs or weak Gibbs
in a rather strong sense (see Remark~\ref{classifyingtool}). In
fact, the importance of the measure~$\mu_*$ lies in the fact that
in a sense it reflects the {\em local} Gibbs structure of $\mu$
(see Proposition~\ref{comparaison}), which proves to be sufficient
for us to show that $\mu$ itself satisfies the multifractal
formalism.
\medskip

\qquad \qquad {\bf 2.4. -- Example: $\beta$-shift of finite type
--} In this section we assume that ${\cal I}_{(\beta,{\bf
d})}=\{\ii_0,\dots,\ii_{{\bf r}-1}\}$, for some ${\bf r}\ge1$. We
also assume that the $\beta$-shift is {\em of finite type}, which
allows us to exhibit a s.a.c naturally associated to $\beta$.
Recall that the $\beta$-shift is of finite type if and only if
there exists $T\ge2$ and $\varepsilon_i\ge0$, such that
$$
1={\varepsilon_1\over\beta^1}+\dots+{\varepsilon_T\over\beta^T}
$$
together with the lexicographic conditions
$$
\varepsilon_i\cdots\varepsilon_{T-1}(\varepsilon_T-1)\varepsilon_1
\cdots\varepsilon_{i-1}
\prec_{\rm lex}\varepsilon_1\cdots\varepsilon_T\qquad (2\le i\le T).
$$
Let ${\bf s}=\sum_{i=1}^T\varepsilon_i$; each $j\in\{\0,\dots,{\bf
s-1}\}$ can be written as follows:
$$
\hbox{$j=\varepsilon_1+\dots+\varepsilon_{k-1}+\varepsilon$ with $1\le k\le
T$ and $0\le
\varepsilon\le\varepsilon_k-1$}.
$$
Put $w_j=\varepsilon_1\dots\varepsilon_{k-1}\varepsilon$. The
s.a.c. ${\cal R}:=\{\widehat{\RR}_j=\RR_{w_j}\}_{j=\0}^{{\bf
s-1}}$ is adapted to [0,1[ because, for any $j=\0,\dots,{\bf
s-1}$, the interval $\widehat{\RR}_j[0,1[$ is bounded from below
by $\displaystyle {\varepsilon_1\over\beta^1}+\cdots+
{\varepsilon_{k-1}\over\beta^{k-1}}+{\varepsilon\over\beta^k}$ and
has length $\beta^{-k}$. For any word $w=\xi_0\dots
\xi_{n-1}\in\{\0,\dots,{\bf s-1}\}^n$, Lemma~\ref{Ath4} yields the
following matrix relation:
\begin{equation}\label{family}
\pmatrix{
\mu(\[w\]+\ii_0)\cr
\vdots\cr
\mu(\[w\]+\ii_{{\bf r}-1})}
=M_{w_{\xi_0}}\cdots M_{w_{\xi_{n-1}}}
\pmatrix{
\mu([0,1[+\ii_0)\cr
\vdots\cr
\mu([0,1[+\ii_{{\bf r}-1})\cr}.
\end{equation}
Consider now the family of the matrices ${\cal
M}:=~\{\widehat{M}_j=M_{w_j}\}_{j=\0}^{\bf s-1}$; according to the
definition of the measures $\mu_k$ in (\ref{defmuj}), the identity
(\ref{family}) turns into
$$
\pmatrix{ \mu_1\[w\]\cr \vdots\cr \mu_{{\bf r}-1}\[w\]}
=\widehat{M}_w \pmatrix{ \mu_1[0,1]\cr \vdots\cr \mu_{{\bf
r}-1}[0,1]\cr},
$$
which means that the $\mu_k$ are indeed ${\cal M}$-measures w.r.t.
the $\cal R$-net.

\end{section}

\begin{section}{
\large
\bf Bernoulli convolution in multinacci bases}

Let $\mu$ be the $(\pp\;\qq)$-distributed $(\beta,2)$-Bernoulli
convolution with $\pp,\qq>0$ and $\beta$ the {\em multinacci
number of degree $\bfm\ge2$}, {\em i.e.}, the PV~number which is
defined as the appropriate root of
$$
\beta^{\bfm}=\beta^{\bfm-1}+\cdots+\beta+1.
$$
The set ${\cal I}_{(\beta,2)}$ from Definition \ref{DefI} has
precisely ${\bf m}+1$ elements, namely, $\ii_0=\0$, $\ii_1=\1$ and
$\ii_k=~\beta^{k-1}-(\beta^{k-2}+\dots+\beta^0)$ for
$k=~2,\dots,\bfm$. Fix $i\in{\cal B}=\{\0,\1\}$; the matrix~$M_i$
defined in (\ref{A9}) is the incidence matrix of the finite
automaton represented in {\it Figure~5}. It has the state set
${\cal I}_{(\beta,2)}$ and the labels $\pp=\pp_0$ and $\qq=\pp_1$.
The arrow from the state $\ii_h$ to the state $\ii_k$ has label
$\pp_j$ if and only if $j=i+\beta\ii_h-\ii_k$ belongs to ${\cal
D}=\{\0,\1\}$. Accordingly,
$$
M_\0=\pmatrix{\pp&0&0&\dots&0\cr0&0&\qq&\dots&0\cr\vdots&\vdots
&\vdots&\ddots&\vdots&\cr
0&0&0&\dots&\qq\cr \qq&\pp&0&\dots&0}\quad\hbox{and}\quad
M_\1=\pmatrix{\qq&\pp&0&\dots&0\cr0&0&0&\dots&0\cr\vdots
&\vdots&\vdots&\ddots&\vdots&\cr
0&0&0&\dots&0\cr 0&\qq&0&\dots&0}.
$$
%
\begin{figure}[!h]
  \begin{center}
       \includegraphics{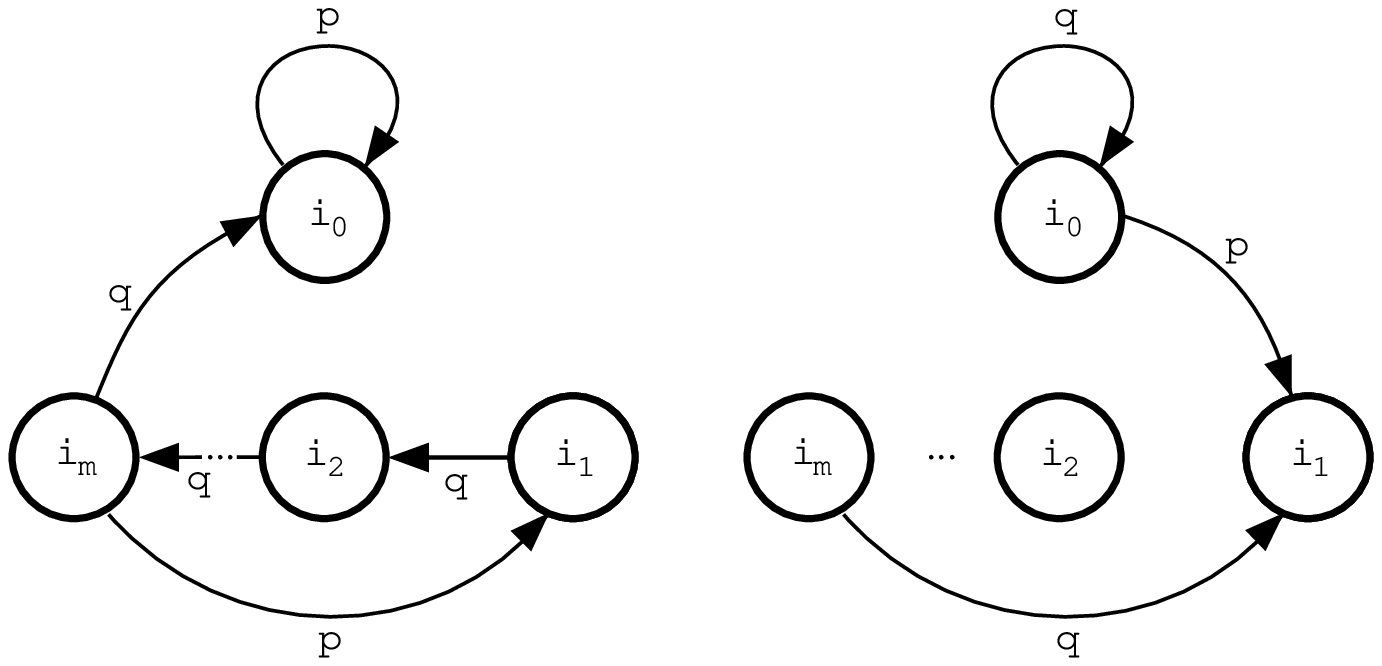}
\caption{\leftskip=1.5cm\rightskip=1.5cm
{ \small\it For $\beta$ being the multinacci number of order ${\bf m}$, we represent the
automaton of the relation
$\cdot
\mathop{\triangleright}\limits^i
\cdot$ on ${\cal I}_{(\beta,2)}$, for $i={\bf 0}$  (left) and $i={\bf 1} (right)$.}}
  \end{center}
\end{figure}
%

\qquad {\bf 3.1. -- The intermediate measure $\mu_*$ --} As will
be shown in Theorem \ref{local3}, there are cases for which $\mu$
cannot be $\gotF$-weak Gibbs for any ``reasonable" $\bf s$-fold
net $\gotF$. However, the multifractal analysis of $\mu$ is still
possible via introducing the intermediate probability measure
$\mu_*$, which turns out to be weak Gibbs (Theorem \ref{thA}) and
equivalent to $\mu$ in a ''strong'' sense (specified in
Proposition \ref{comparaison}). The measure $\mu_*$ defined in
(\ref{defPiprime}) below, will be given by an application of the
process described in section 2.3.


\qquad To begin with, notice that the $\beta$-shift associated to
the multinacci number is of finite type, whence, by the result of
Section~2.4, one can find a finite family of words $\{w_j\}_j$ in
the alphabet $\{\0,\1\}$ for which the system of affine
contractions $\{\widehat{\RR}_j=\RR_{w_j}\}_j$ is adapted to the
unit interval. However, for the sake of simplicity, our family of
words will be different from the one obtained by the systematic
approach presented in Section~2.4. This way we obtain a family of
just $2\times2$ matrices.

\qquad Heuristically, our approach is based on the following two
remarks. Denoting by
$$
X=\pmatrix{1&0&0&0&\!\!\!\cdots\!\!\!&0\cr0&0&1&0&\!\!\!\cdots\!\!\!&0\cr},
$$
one obtain first
$$
XM_{\0^{\bfm}}=\pmatrix{{\pp}^{\bfm}&0\cr\pp{\qq}^{\bfm-1}&\pp{\qq}^{\bfm-1}}X\;.
$$
Furthermore, since\footnote{Given two square matrices
$A=(a_{i,j})$ and $B=(b_{i,j})$ with nonnegative entries, we write
$A\prec B$ whenever $b_{i,j}=0$ implies that $a_{i,j}=0$, for any
$i,j$.} $ M_{\0^q\1^r\1}\prec \pmatrix{
1&1\!\!&0&\!\!\!\!\cdots\!\!\!\!&0\cr
\vdots&\vdots\!\!&\vdots&\!\!\!\!\ddots\!\!\!\!&\vdots\cr
1&1\!\!&0&\!\!\!\!\cdots\!\!\!\!&0\cr} $, for any $q,r\ge0$, it is
also clear that
$$
XM_{\0^q\1^{r}\1\0}\prec
\pmatrix{1&1\cr1&1}X.
$$
Actually, there exists a family of words $w_j$ of the form either
$\0^{\bfm}$ or $\0^x\1^y\1\0$, for which the system
$\{\RR_{w_j}\}_j$ is adapted to the unit interval. To see this,
notice that the algebraic property  of $\beta$ ensures that the
semi-open intervals $\RR_{\1^y\1\0}[0,1[$, for $0\le y\le \bfm-2$,
form a partition of $[1/\beta,1[$. Moreover,
$\RR_{\0^{\bfm}}\left[0,1\right[$ and
$\RR_{\0^x}\left[1/\beta,1\right[$, for $0\le x\le \bfm-1$, form a
partition of $[0,1[$. For a suitable labelling of the words
$\0^x\1^y\1\0$,  let $j=j(x,y)$ be the integer such that
\begin{equation}
j-1:=(\bfm-1)\Big((\bfm-1)-x\Big)+y\;.\label{BB2}
\end{equation}
The uniqueness of the division of $j-1$ by $\bfm-1$ with remainder
implies $(x,y)\mapsto j$ is a bijection from
$\big\{0,\dots,\bfm-1\big\}\times\big\{0,\dots,\bfm-2\big\}$ onto
$\big\{1,\dots,\bfm(\bfm-1)\big\}$. Next, we consider the adapted
system ${\cal
R}:=\{\widehat{\RR}_j=\RR_{w_j}\}_{j=0}^{\bfm(\bfm-1)}$, where
$$
w_0:=\0^{\bfm}\;\;\hbox{and}\;\;w_{j}:=\0^{\xx(j)}\
\1^{\yy(j)}\1\0\quad\hbox{with}\quad
\cases{
\xx(j):=\displaystyle(\bfm-1)-\left[{j-1\over \bfm-1}\right]\cr\cr
\yy(j):=\displaystyle(\bfm-1)\left\{{j-1\over \bfm-1}\right\}}
$$
(here,  $[\;\cdot\;]$ and $\{\;\cdot\;\}$ stand respectively for
the integral and the fractional~part of a number).

\qquad Put
$Y=\pmatrix{1&0&0&\!\!\!\cdots\!\!\!&0\cr0&1&0&\!\!\!\cdots\!\!\!&0\cr}$;
given two integers $x$ and $y$, one has successively:
$$
XM_{\0^x\1}=\cases{
\pmatrix{\pp^x\qq&\pp^{x+1}\cr0&0}Y & if $x<\bfm-2$;\cr\cr
\pmatrix{\pp^{\bfm-2}\qq&\pp^{\bfm-1}\cr0&\qq^{\bfm-1}}Y & if $x=\bfm-2$;\cr\cr
\pmatrix{\pp^{\bfm-1}\qq&\pp^{\bfm}\cr\qq^\bfm&\qq^{\bfm-1}\pp}Y & if
$x=\bfm-1$;\cr}
$$
and
$$
YM_{\1^y\0}=\cases{
\pmatrix{\pp&0\cr0&\qq}X & if $y=0$;\cr\cr
\pmatrix{\qq^y\pp&\qq^y\pp\cr0&0}X & if $y>0$.\cr}\qquad\qquad\,
$$
Put $\alpha=(\qq/\pp)^{\bfm -1}$; a straightforward computation
yields the following
\begin{lemma}\label{defPi}
{\sl For any $j\in{\cal J}:=\Big\{\0,\dots,\bfm(\bfm-1)\Big\}$ one
has $XM_{w_j}=P_jX$, where
$$
\matrix{
\bullet\; P_\0:=\pp^\bfm\pmatrix{1&0\cr\alpha&\alpha}
\hfill&
\hbox{when $j=\0$ ({\it i.e.,}
$w_{\bf j}=\0^\bfm$);}
\hfill\cr\cr
\bullet\;P_j:=\pp^\bfm\qq^j\pmatrix{1&
1\cr\alpha&\alpha}
\hfill&
\hbox{when $\0< j<\bfm$ ({\it i.e.,}
$w_j=\0^{\bfm-1}\1^{j-1}\1\0$);}
\hfill\cr\cr
\bullet\;P_\bfm:=\qq^\bfm\pmatrix{1/\alpha
&1/\alpha\cr0&1}
\hfill&
\hbox{when $j=\bfm$ ({\it i.e.,} $w_j=\0^{\bfm-2}\1\0$);}
\hfill\cr\cr
\bullet\;P_j:=\pp^{\xx(j)+1}\qq^{\yy(j)+1}\pmatrix{1&1\cr0&0}
\hfill&
\hbox{when $\bfm<j\le\bfm(\bfm-1)$.}
\hfill}
$$
}
\end{lemma}

From here on, we will consider the basic intervals associated to
the adapted system ${\cal R}$. For any word $w\in{\cal J}^*$ we
denote $\[w\]=\widehat{\RR}_w[0;1[$; by (\ref{correct2}) and
Lemma~\ref{defPi},
$$
\pmatrix{\mu_0\[w\]\cr\mu_2\[w\]}=
P_w\pmatrix{\mu_0([0,1])\cr\mu_2([0,1])}.
$$
(Recall that $\mu_i(\, \cdot \,):=\mu(\,\cdot\cap[0,1]+\ii_i)$.)
We introduce the probability
\begin{equation}\label{defPiprime}
\mu_*:={1\over \mu_0([0,1])+\mu_2([0,1])}\,(\mu_0+\mu_2),
\end{equation}
so that, for any  word $w\in{\cal J}^*$,
$$
\mu_*\[w\]=LP_wR,
$$
where
$$
L=\pmatrix{1&1}\quad\hbox{and}\quad R={1\over
\mu_0([0,1])+\mu_2([0,1])}\pmatrix{\mu_0([0,1])\cr\mu_2([0,1])}
=\pmatrix{1-\qq^{\bfm-1}\cr\qq^{\bfm-1}}.
$$
Now we are ready to formulate one of the central claims of the
present paper:

\begin{theorem}\label{thA}
{\sl The Gibbs properties of $\mu_*$ w.r.t. the net $\gotF$
associated to the adapted system $\cal R$, are the following:
\hfil\break \null\qquad (i) when ${\bf m}=2$, the measure $\mu_*$
is weak Gibbs if $\pp=\qq$ and Gibbs if $\pp\ne\qq$; \hfil\break
\null\qquad (ii) when ${\bf m}\ge3$, the measure $\mu_*$ is weak
Gibbs if $\pp=\qq$ and Gibbs if $\pp>\qq$.}
\end{theorem}

\begin{remark}\label{truc}{The measure $\mu_*$ is not even
$\gotF$-weak Gibbs when ${\bf m}\ge3$ and $\pp<\qq$: consider the
basic interval $\[\xi_0\dots\xi_{2n-1}\]=\[{\bf m}^n({\bf
m}+1)^n\]$; then, it is clear that the ratio
$$
{\mu_*\[\xi_0\dots\xi_{2n-1}\]\over\mu_*\[\xi_0\dots\xi_{n-1}\]
\mu_*\[\xi_n\dots\xi_{2n-1}\]}
$$
does not satisfy the condition
$$
{1\over K_n}\le{\mu_*\[\xi_0\dots\xi_{2n-1}\]
\over\mu_*\[\xi_0\dots\xi_{n-1}\]\mu_*\[\xi_n\dots\xi_{2n-1}\]}
\le K_n
$$
together with $\lim_n{1\over n}\log K_n=0$ (this will be explained
in more detail in Remark~\ref{discontinuous}). In the case ${\bf
m}\ge3$ and $\pp<\qq$, the problem of the existence of a net
$\gotF'$ of $[0,1]$, w.r.t. which $\mu_*$ is $\gotF'$-weak Gibbs,
remains open.}
\end{remark}

Theorem~Ê\ref{thA} is a corollary of Theorem~\ref{Holder} below,
whose formulation needs the explicit expression for a potential
$\Phi$ associated with~$\mu_*$.

\qquad{\bf 3.2. -- The potential associated to $\mu_*$ --} In
order to compute the limit potential of the $n$-step potential
associated to the measure $\mu_*$, we need some classical facts
about continued fractions which can be found in \cite{Per1}. Let
$u_0,u_1,\dots$ and $v_0,v_1,\dots$ be two infinite sequences of
real numbers, both assumed to be positive, except, possibly, $u_0$
which is allowed to be nonnegative. For any $n\ge0$ the positive
reals $p_n$ and $q_n$ are defined by
\begin{equation}\label{Khint1}
\pmatrix{p_n\cr
q_n}=\pmatrix{u_0&v_0\cr1&0}\pmatrix{u_1&v_1\cr1&0}
\cdots\pmatrix{u_n&v_n\cr1&0}\pmatrix{1\cr 0},
\end{equation}
and by convention, $\pmatrix{p_{-1}\cr q_{-1}}=\pmatrix{1\cr 0}$.
Then $\displaystyle{p_n/q_n}$ is the {\em continued fraction}
associated with $u_0,\dots,u_n$ and $v_0,\dots,v_n$ in the usual
sense: for any $n\ge0$,
\begin{equation}\label{Khint2}
{p_n\over q_n}=u_0+
{{v_0}\over{\displaystyle u_1+
{{v_1\raise10pt\hbox{}}\over{\displaystyle\;\;\,\raise-8pt\hbox{$\ddots$}\;\;\,
\raise-20pt\hbox{$\displaystyle
u_{n-1}+ {v_{n-1}\over u_n}$}
\!\!\!\!\!\!\!\!\!\!\!\!}}}
\!\!\!\!\!\!\!\!\!}\qquad\qquad
\end{equation}
(see \cite{Per1} for the proof). A direct consequence of
(\ref{Khint1}) is the following relation:
\begin{equation}\label{Khint3}
\pmatrix{u_0&v_0\cr1&0}\cdots\pmatrix{u_n&v_n\cr1&0}=
\pmatrix{p_n&v_np_{n-1}\cr q_n&v_nq_{n-1}}.
\end{equation}
Also, for $n\ge1$, we have
\begin{equation}\label{Khint5}
\cases{
p_n=u_np_{n-1}+v_{n-1}p_{n-2}\cr
q_n=u_nq_{n-1}+v_{n-1}q_{n-2}\ .\cr}
\end{equation}

%

\qquad Consider a special case of such continued fractions. Recall that
$\alpha=(\qq/\pp)^{\bfm -1}$; given a sequence of integers
$a_0,a_1,\dots$ with $a_0\ge0$ and $a_i>0$ for $i\ge 1$ and
$\kappa=\0$ or $\1$, put
\begin{equation}\label{product1}
\cases{u_n=\alpha^1+\dots+\alpha^{a_n}\qquad\;\;\hbox{and}\quad
v_n=\alpha^{a_n}&if $n+\kappa$ is even\cr\cr
u_n={1/\alpha^1}+\dots+{1/\alpha^{a_n}}\;\;\hbox{and}\quad
v_n={1/\alpha^{a_n}} & if $n+\kappa$ is odd.}
\end{equation}
By convention, put $u_0=0$ in either case when $a_0=0$. Then the
sequences $(p_n)_{n=-1}^\infty$ and $(q_n)_{n=-1}^\infty$ are
defined by (\ref{Khint5}). The sequence
$(p_{2n}/q_{2n})_{n=0}^\infty$ is nondecreasing while
$(p_{2n+1}/q_{2n+1})_{n=0}^\infty$  is nonincreasing, and
$p_{2n}/q_{2n}\le p_{2n+1}/q_{2n+1}$. Moreover, setting
$\rho:=\min\{\sqrt{\alpha},1/\sqrt{\alpha}\}$, for any $n\ge1$,
one has
\begin{equation}\label{convcontfrac}
\left\vert{p_n\over q_n}-{p_{n-1}\over q_{n-1}}\right\vert\le
\cases{\displaystyle{\rho^{a_1+\dots+a_n}\over \rho^{2+2a_0+a_n}}&
if $\alpha\ne1$\cr \displaystyle{1\over a_0+\dots+a_n}& if
$\alpha=1$.}
\end{equation}
This is a property of the regular continued fractions when
$\alpha=1$ (see, e.g., \cite{Kin}) and corresponds to part~(ii) of
Lemma~\ref{Holder2} in case $\alpha\ne1$---see below. Analogously
to the regular continued fractions, we denote
$$
\left[\kappa\vert a_0;a_1,\dots,a_n\right]={p_n\over
q_n}\quad\hbox{and}\quad\cases{\left[\kappa\vert
a_0;a_1,\dots,a_n,\infty\right]=
\displaystyle\lim_{k\to\infty}\left[\kappa\vert
a_0;a_1,\dots,a_n,k\right]\cr\cr\left[\kappa\vert
a_0;a_1,a_2,\dots\right]= \displaystyle\lim_{n\to\infty}
\left[\kappa\vert a_0;a_1,\dots,a_n\right].}
$$
We also define the matrix
$$
Q_\kappa(a_0,\dots,a_{n})
=\pmatrix{u_0&v_0\cr1&0}\cdots\pmatrix{u_n&v_n\cr1&0}
$$
and, for any column vector $\pmatrix{x\cr y}$ with nonnegative
entries,
\begin{equation}\label{xy}
\left[\kappa\vert a_0;a_1,\dots,a_n\vert \pmatrix{x\cr
y}\right]={p_n(x,y)\over
q_n(x,y)}\quad\hbox{with}\quad\pmatrix{p_n(x,y)\cr q_n(x,y)}=
Q_\kappa(a_0,\dots,a_{n})\pmatrix{x\cr y}.
\end{equation}

Notice that if $ \Delta$ denotes $\pmatrix{0&1\cr1&0}$, then for
any $n\ge0$,
$$
\pmatrix{
u_n&v_n\cr1&0\cr}=\cases{\displaystyle\Delta\left({1\over\pp^{\bf m}} P_\0\right)^{a_n},&if $n+\kappa$ is
even and\cr\cr
\displaystyle\left({1\over \qq^{\bf m}} P_{\bf m}\right)^{a_n}\!\!\Delta,&if $n+\kappa$ is odd.}
$$

This allows to compute the potential $\Phi:{\cal J}^{\bf N}\to{\bf
R}$ associated with $\mu_*$ . To do so, let us first introduce a
suitable notation for the words in $\{\0,\bfm\}^*$: given $\xi=\0$
or $\bfm$, we denote by $\osc{\xi}{a}$ the word $\xi^a$, for any
nonnegative integer $a$. Furthermore, for any sequence of
nonnegative integers $a_1,\dots,a_n$ ($n\ge2$), we define the word
$\osc{\xi}{a_1,\dots,a_n}$ by the induction relation
$$
\osc{\xi}{a_1,\dots,a_n}=\osc{\xi}{a_1}\osc{\bfm-\xi}{a_2,\dots,a_n}.
$$
Before proving the uniform convergence of the $n$-step potential
associated with $\mu_*$, we give an explicit formula for the limit
potential which we denote by $\Phi$. Given $j\in{\cal J}$, we
distinguish between the same cases as in Lemma~\ref{defPi}. Put

\null\qquad$\bullet$ $\;\displaystyle
\Xi(j):=\cases{\pp^\bfm\alpha^2 & if $j=\0$ ;\cr
\pp^\bfm\qq^j(1+\alpha)& if $\0< j<\bfm$ ;\cr
\qq^{\bfm}/\alpha^2& if $j=\bfm$ ;\cr
\pp^{\xx(j)+1}\qq^{\yy(j)+1}
&  if $\bfm<j\le\bfm(\bfm-1)$\cr}$

\null and when $0<j<\bfm$,

\null\qquad$\bullet$ $\;\displaystyle
X_j:=\cases{\pmatrix{1\cr\alpha}&
if $\0<j<\bfm\;;$\cr
\pmatrix{1\cr0}& if $\bfm<j\le \bfm(\bfm-1)$.\cr}$

\qquad Given $\kappa\in\{\0,\1\}$, we denote
$\widehat\kappa=1-\kappa$ and $\kappa\star\bfm$ stands for either
$\0$ or $\bfm$ when $\kappa$ is either $\0$ or $\1$ respectively;
hence, $\kappa\star\bfm$ is just "multiplication" of $\kappa$ by
$\bfm$, which is not be confused with the concatenation
$\kappa\bfm$.

\qquad Let ${\cal J}^{\bf N}\!\ni\!\xi\!=\!wj\omega$, with
$\{\0,\bfm\}^*\!\ni\! w\!=\!\osc{\kappa\star
\bfm}{a_1,\dots,a_n}$, $j\in{\cal J}\backslash\{\0,\bfm\}$ and
$\omega\!\in\!{\cal J}^{\bf N}$. Assume $\kappa=\0$ (the case
$\kappa=\1$ is symmetric); for $k\ge1+a_1+\dots+a_n$ direct
computation yields
\begin{eqnarray*}
\Phi(\xi)=\phi_k(wj\omega)&=&
\log(\pp^\bfm)+\log\left({\pmatrix{1&1}
Q_\0(a_1,\dots,a_n)
\Delta^{1+n}X_j\over
\pmatrix{1&1}Q_\0(a_1-1,\dots,a_n)\Delta^{1+n}X_j}\right)\\&=&
\log\Xi(\0)+\log\left({\pmatrix{1&0}
Q_\1(1,a_1,\dots,a_n)
\Delta^{1+n}X_j\over
\pmatrix{0&1}Q_\1(1,a_1,\dots,a_n)\Delta^{1+n}X_j}\right).
\end{eqnarray*}
More generally, for $\kappa\in\{\0,\1\}$, $j\in{\cal
J}\backslash\{\0,\bfm\}$ and $\omega\in{\cal J}^{\bf N}$,
\begin{eqnarray*}
&\bullet&\!\Phi(wj\omega)\!=\!\displaystyle\log\Big({\Xi(\kappa\star\bfm)
\left[\widehat\kappa\vert1;a_1,\dots,a_n\vert
\Delta^{\widehat\kappa+n}X_i\right]}\Big)\quad\hbox{if }
w=\osc{\kappa\star\bfm}{a_1,\dots,a_n};\\\\
&\bullet&\Phi(j\omega)\;=\;\log\Big({\Xi(j)}\Big);
\\\\
&\bullet&\Phi(\omega)\;=\;\cases{
\displaystyle\log\Big({\Xi({\kappa\star\bfm})
[\widehat\kappa\vert1;a_1,\dots,a_n,\infty]}\Big) & if
$\omega=\osc{\kappa\star\bfm}{a_1,\dots,a_n,\infty}$\cr\cr
\displaystyle \log\Big({\Xi({\kappa\star\bfm})
[\widehat\kappa\vert1;a_1,a_2,\dots]}\Big) & if
$\omega=\osc{\kappa\star\bfm}{a_1,a_2,\dots}$.\cr}
\end{eqnarray*}

\qquad We have sketched the proof of the pointwise convergence of
the $n$-step potential $\phi_k$ to the potential $\Phi$ whose
expression is given above. Indeed, by (\ref{convcontfrac}), this
is a simple consequence of the convergence of the continued
fractions involved. We need however to show that the convergence
is uniform; this is dealt with in the following theorem, by means
of considering the two alternative cases \hbox{$\pp\ne\qq$ and
$\pp=\qq$.}

\begin{theorem}\label{Holder}
{\sl Let $\rho:=\min\{\sqrt\alpha,1/\sqrt\alpha\}\le1$; there
exists a constant $K>0$ such that, for any $\omega\in{\cal J}^{\bf
N}$ and any $n\ge1$, we have
\begin{equation}\label{expconv}
\Vert\phi_n-\Phi\Vert_\infty\le \cases{ K\rho^n&if $\pp>\qq$ \cr
K\rho^n&if $\pp\ne\qq$ and $\bfm=2$\cr {K/ n}&if $\pp=\qq$.}
\end{equation}
}
\end{theorem}

\begin{remark}\label{discontinuous}{The potential $\Phi$ is
continuous except when $\pp<\qq$ and $\bfm\ge3$, as in that case
$\Phi$ is discontinuous at $\omega={\bfm}^\infty$. Indeed, as we
have already noticed (see Remark~\ref{truc}), $\mu_*$ is not
$\gotF$-weak Gibbs if $\pp<\qq$ and $\bfm\ge3$. Let us also point
out that $\Phi$ is H\"older continuous if either $\pp>\qq$ or
$\pp\ne\qq$ and $\bfm=2$, yielding $\mu_*$ to be $\gotF$-Gibbs in
these cases.}
\end{remark}

\qquad {\bf 3.3. --  Preliminary to proof of Theorem~\ref{Holder}
--} The key argument is based on the fact that the pointwise
convergence of the $n$-step potential to a continuous limit
implies the uniform convergence. To see this, let us define the
$n$-step variation of an arbitrary map $f:\Omega:={\cal A}^{\bf
N}\to{\bf R}$ with ${\cal A}:=\{\0,\dots,{\bf s-1}\}$, by defining
$$
\hbox{Var}_n(f)=\sup\Big\{f(\xi)-f(\xi^0)\;;\;(\xi,\xi^0)\in{\cal
A}^{\bf N}\times{\cal A}^{\bf
N}\;\hbox{and}\;\xi_0\dots\xi_{n-1}=\xi^0_0\dots\xi^0_{n-1}\Big\}.
$$
It is clear that $f$ is continuous if and only if
$\lim_{n\to\infty}\hbox{Var}_n(f)=0$.

\begin{lemma}\label{pointwise}
{\sl Suppose that  ${\cal S}$ is an s.a.c. adapted to the interval
$[0,1]$ and let $\eta$ be a probability measure whose support is a
subset of $[0,1]$. If the $n$-step potential $\phi_n$ associated
to $\eta$ converges to $\phi:\Omega\to{\bf R}$ pointwise, then
$\Vert \phi_n-\phi\Vert_\infty\le\hbox{Var}_n(\phi)$.}
\end{lemma}
{\bf Proof.} Given an arbitrary rank $n$ and any $\varepsilon>0$,
the pointwise convergence of $\phi_n$ to $\phi$ implies that for
any $\xi\in\Omega$, there exists an integer $N(\xi)\ge n$ such
that
\begin{equation}\label{pointw}
|\phi_{N(\xi)}(\xi)-\phi(\xi)|\le\varepsilon.
\end{equation}
Since the product space $\Omega$ is compact, there exists a finite
set $X\subset\Omega$ such that
\begin{equation}\label{union}
\Omega=\bigcup_{\xi\in X}\[\xi_0\dots\xi_{N(\xi)-1}\].
\end{equation}
In the product space $\Omega$, the intersection of any pair of
cylinders is either empty or coincides with one of them, whence we
may take a set smaller than $X$ so that the union in (\ref{union})
becomes disjoint. For any $\omega\in\Omega$ put
$X_{\omega,n}=X\cap[\omega_0\dots\omega_{n-1}]$. Then
$$
\exp\left(\phi_n(\omega)\right)={\eta\[\omega_0\dots\omega_{n-1}\]
\over \eta\[\omega_1\dots\omega_{n-1}\]}={\sum_{\xi\in
X_{\omega,n}} \eta\[\xi_0\dots\xi_{N(\xi)-1}\]\over \sum_{\xi\in
X_{\omega,n}}\eta\[\xi_1\dots\xi_{N(\xi)-1}\]},
$$
whence
$$
\min_{\xi\in X_{\omega,n}}\left\{{\eta\[\xi_0\dots\xi_{N(\xi)-1}\]\over
\eta\[\xi_1\dots\xi_{N(\xi)-1}\]}\right\}\le
\exp\left(\phi_n(\omega)\right)\le\max_{\xi\in
X_{\omega,n}}\left\{{\eta\[\xi_0\dots\xi_{N(\xi)-1}\]\over
\eta\[\xi_1\dots\xi_{N(\xi)-1}\]}\right\}.
$$
Since by definition
$$
\displaystyle{\eta\[\xi_0\dots\xi_{N(\xi)-1}\]\over
\eta\[\xi_1\dots\xi_{N(\xi)-1}\]}
=\exp\left(\phi_{N(\xi)}(\xi)\right),
$$
we obtain, in view of (\ref{pointw}),
$$
\min_{\xi\in X_{\omega,n}}\left\{\phi(\xi)\right\}-\varepsilon\le
\phi_n(\omega)\le\max_{\xi\in
X_{\omega,n}}\left\{\phi(\xi)\right\}+\varepsilon
$$
and thus,
\begin{equation}\label{osci}
\phi(\omega)-\hbox{Var}_n(\phi)-\varepsilon\le
\phi_n(\omega)\le\phi(\omega)+\hbox{Var}_n(\phi)+\varepsilon.
\end{equation}
Since (\ref{osci}) holds for an arbitrary $\varepsilon>0$, we have
$\vert\phi_n(\omega)-\phi(\omega)\vert\le\hbox{Var}_n(\phi)$,
which concludes the proof.\hfill\break \null\hfill\qed
\medskip

\qquad In what follows, the sequences $(u_n)_{n=0}^\infty$,
$(v_n)_{n=0}^\infty$, $\dots$  are always associated with
$\kappa\in\{\0,\1\}$ and $a_0,a_1,\dots$, as in
(\ref{product1}), (\ref{Khint1}) and (\ref{xy}). For any $n\ge1$,
we introduce the following quantity:
$$
\delta_n:=\left\vert{p_n\over q_n}-{p_{n-1}\over
q_{n-1}}\right\vert.
$$
We need two lemmas which involve $\delta_n$.

\begin{lemma}\label{Holder1}
{\sl For any integer $n\ge1$ and any column vector with
nonnegative entries $\pmatrix{x\cr y}\neq\pmatrix{0\cr 0}$, one
has:
$$
\hbox{(i) : }\;\left\vert{p_n(x,y)\over q_n(x,y)}-{p_n\over q_n}\right\vert
\le{yv_n\over xu_n+yv_n}\cdot \delta_n\quad\hbox{and}\quad\hbox{(ii) : }\;\delta_{n+1}
\le{v_n\over u_nu_{n+1}+v_n}\cdot \delta_n.
$$}
\end{lemma}

{\bf Proof.} (i) By (\ref{Khint3}),
\begin{eqnarray*}
{p_n(x,y)\over q_n(x,y)}
\!\!\!&=&\!\!\!{\pmatrix{1&0}\pmatrix{p_n&v_np_{n-1}\cr
q_n&v_nq_{n-1}}
\pmatrix{x\cr y}\over
\pmatrix{0&1}\pmatrix{p_n&v_np_{n-1}\cr
q_n&v_nq_{n-1}}\pmatrix{x\cr y}}\\
\!\!\!&=&\!\!\!{xp_n+yv_np_{n-1}\over
xq_n+yv_nq_{n-1}}\\
\!\!\!&=&\!\!\!{aq_n\over xq_n+y{v_nq_{n-1}}}\cdot {p_n\over
q_n}+{y{v_nq_{n-1}}\over xq_n+ y{v_nq_{n-1}}}\cdot{p_{n-1}\over
q_{n-1}},
\end{eqnarray*}
whence
\begin{eqnarray*}
{p_n(x,y)\over q_n(x,y)}-{p_n\over q_n}
\!\!\!&=&\!\!\!{y{v_nq_{n-1}}\over xq_n+
y{v_nq_{n-1}}}\cdot\left({p_{n-1}\over q_{n-1}}-{p_n\over q_n}\right)
\end{eqnarray*}
and since $q_n\ge u_nq_{n-1}$, we are done.

(ii) We obtain (ii) by simply applying (i) to the vector
$\pmatrix{x\cr y}=\pmatrix{u_{n+1}\cr 1}$.\hfill\break
\null\hfill\qed

\begin{lemma}\label{Holder2}
{\sl (i) : For an arbitrary rank $n\ge1$ we have
$\delta_{n+1}\le\delta_n$, and
$$
\delta_n\le{\rho^{a_1+\dots+a_{n}}\over \rho^{2+2a_0}}
$$
in either of the following cases: $n+\kappa$ is even and
$\alpha>1$; or $n+\kappa$ is odd and $\alpha<1$;\hfil\break
\null\qquad (ii) : if $\alpha\ne1$, then for any $n\ge1$,
$$
\delta_n\le{\rho^{a_1+\dots+a_{n}}\over \rho^{2+2a_0+a_n}}\;;
$$
\null\qquad (iii) : if $\alpha\ne1$, then there exists a constant
$K>0$ such that, for arbitrary rank $n\ge1$ and any
integer~$0<a<a_n$,
$$
\Big\vert [\kappa\vert a_0;a_1,\dots,a_{n-1},a_n]
-[\kappa\vert a_0;a_1,\dots,a_{n-1},a]\Big\vert\le
{K\over \rho^{2a_0}}\,
\rho^{a_1+\dots+a_{n-1}+a}.
$$
}
\end{lemma}

{\bf Proof.} (i) : We are going to establish the inequality in (i)
in the case when $n+\kappa$ even with $\alpha>1$; this implies the
inequality in the opposite case, as the value of $\delta_n$
remains the same whenever a pair $(\alpha, \kappa)$ is replaced by
$(\alpha'=1/\alpha>1, \widehat\kappa=1-\kappa)$.

Assume now $n+\kappa$ is even, and let $k\le n$ be such that
$k+\kappa$ is even as well. Then, on one hand, Lemma
\ref{Holder1}~(ii) implies
$$
{\delta_k\over \delta_{k-1}}\le{v_{k-1}\over u_{k-1}u_k}
\le
{1/\alpha^{a_{k-1}}\over
(1/\alpha)\alpha^{a_k}}
\le{1\over \alpha^{(a_{k}+a_{k-1})/2}}\;,
$$

and on the other hand, since $\displaystyle{\delta_{k-1}/
\delta_{k-2}}\le1$,
$$
\delta_{n}\le{\delta_n\over\delta_{n-1}}\,
{\delta_{n-2}\over\delta_{n-3}}
\cdots{\delta_{\kappa+2}\over\delta_{\kappa+1}}\,\delta_1
\le{\delta_1\over \alpha^{(a_n+a_{n-1}+\dots+a_{\kappa+1})/2}}.
$$
This proves (i), as $\delta_1={v_0/ u_1}$ is bounded by
$\alpha^{1+a_0}$ if $\kappa=0$ and by
$1/\alpha^{a_0+a_1}\le~\alpha^{1+a_0}/\alpha^{a_1/2}$ if
$\kappa=1$.

\qquad (ii) and (iii) : Part (ii) is a straightforward consequence
of (i) and the definition of $\rho$. In order to prove (iii), we
consider, for any integer $0<a<a'$, the quantity
$$
\delta_n(a,a'):=\Big\vert [\kappa\vert a_0;a_1,\dots,a_{n-1},a']
-[\kappa\vert a_0;a_1,\dots,a_{n-1},a]\Big\vert
\le\sum_{i=a}^{a'-1}\delta_n(i,i+1).
$$
Since $[\kappa\vert a_0;a_1,\dots,a_{n-1},i+1] =[\kappa\vert
a_0;a_1,\dots,a_{n-1},i,1]$, we obtain, in view of (ii),
$$
\delta_n(i,i+1)
\le{\rho^{a_1+\dots+a_{n-1}+i}/ \rho^{2+2a_0}},
$$
whence
$$
\delta_n(a,a_n)\le\sum_{i=a}^{a_n-1}\delta_n(i,i+1)\le
\left\{\sum_{j=0}^\infty\rho^{j-2}\right\}{1\over
\rho^{2a_0}}\,\rho^{a_1+\dots+a_{n-1}+a}.
$$
The inequality in (iii) follows from the fact that
$\rho<1$.\hfill\break \null\hfill\qed
\medskip

{\bf Proof of Theorem~\ref{Holder}.} In view of
Lemma~\ref{pointwise}, it suffices to establish the desired
estimate of  $\hbox{Var}_n(\Phi)$. Fix $w\in{\cal J}^n$;  if
$w\notin\{\0,{\bf m}\}^n$, then
$$
\sup\Big\{\vert\Phi(\xi)-\Phi(\xi')\vert\;;\;\xi,\xi'\in\[w\]\Big\}=
~0,
$$
so from here on we assume that $w\in\{\0,{\bf m}\}^n$. Let
$\xi\in\[w\]$ be of the form $\xi=w'j\omega$, where
$$
w'=\osc{\kappa\star\bfm}{a_1,\dots,a_k}
$$
(with $\kappa=\0$ or $\1$), $j\in{\cal J}\backslash\{\0,\bfm\}$
and $\omega\in{\cal J}^{\bf N}$ (the case of $\xi\in\{\0,{\bf
m}\}^{\bf N}$ is handled in the same way). Then there exist
$0<k_n\le k$ and $0< a_{k_n}'\le a_{k_n}$ such that
$$
w=\xi_0\dots\xi_{n-1}=
\osc{\kappa\star\bfm}{a_1,\dots,a_{k_n-1},a_{k_n}'}
$$
(note that $n=a_1+\cdots+a_{k_n}'$). Put
$$
\phi^w=\log\Big({\Xi({\kappa\star\bfm
})[\widehat\kappa\vert1;a_1,\dots,a_{k_n-1},a_{k_n}']}\Big).
$$
We have
\begin{equation}\label{differences}
\Phi(\xi)-\phi^w=A_k+B_k,\quad\hbox{where}\quad
\cases{\displaystyle A_k:=
\log\left({[\widehat\kappa\vert1;a_1,\dots,a_k\vert\Delta^{\widehat\kappa+k}X_j]\over
[\widehat\kappa\vert1;a_1,\dots,a_{k_n-1},a_{k_n}]}\right),\cr
\displaystyle
B_k:=\log\left({[\widehat\kappa\vert1;a_1,\dots,a_{k_n-1},a_{k_n}]\over
[\widehat\kappa\vert1;a_1,\dots,a_{k_n-1},a_{k_n}']}\right)\cr}
\end{equation}
(notice that $\phi^w=\Phi(w({\bf m}+1)w({\bf m}+1)\dots)$ in the
case ${\bf m}\ge 3$). Lemma~\ref{Holder2}~(iii) now yields
\begin{equation}\label{comp1}
B_k\leq
K\,\rho^{a_1+\cdots+a_{k_n-1}+a_{k_n}'}/\rho^2=K\,\rho^{n-2}.
\end{equation}
In order to establish a suitable upper bound of $|A_k|$ in
(\ref{differences}), we first note that if $k>k_n$, then by
definition,
$$
[\widehat\kappa\vert1;a_1,\dots,a_k\vert\Delta^{\widehat\kappa+k}X_j]=
[\widehat\kappa\vert1;a_1,\dots,a_{k_n},a_{k_n+1}\vert Y],
$$
for some nonnegative column vector $Y$. Therefore,
$[\widehat\kappa\vert1;a_1,\dots,a_k\vert\Delta^{\widehat\kappa+k}X_j]$
lies between $[\widehat\kappa\vert1;a_1,\dots,a_{k_n}]$ and
$[\widehat\kappa\vert1;a_1,\dots,a_{k_n},a_{k_n+1}]$ and
Lemma~\ref{Holder2}~(ii) yields
\begin{equation}\label{comp2}
|A_k|\leq \rho^{a_1+\cdots+a_{k_n}}/\rho^4\le\rho^{n-4}.
\end{equation}
It remains to establish the upper bound for $|A_k|$ in the case
$k=k_n$. Note that by Lemma~\ref{Holder1},
\begin{equation}\label{delta_k}
|A_k|\leq\rho^2\Big\vert[\widehat\kappa\vert1;a_1,
\dots,a_k\vert\Delta^{\widehat\kappa+k}X_j]
-[\widehat\kappa\vert1;a_1,\dots,a_k]\Big\vert
\le\rho^2\cdot{v_kx_2\over u_kx_1+v_kx_2}\cdot \delta_k,
\end{equation}
where $x_1$ and $x_2$ are the coordinates of
$\Delta^{\widehat\kappa+k}X_j$. Consider three different cases.

\quad$\bullet$ {\it The case $\alpha<1$. -- } On one hand, if
$\widehat\kappa+k$ is odd, then
\begin{equation}\label{comp3}
|A_k|\le\rho^2\delta_k\le \rho^{a_1+\cdots+a_k}/\rho^2.
\end{equation}
On the other hand, if $\widehat\kappa+k$ is even, then there
exists a constant $K'>0$ such that
\begin{equation}\label{comp4}
|A_k|\le
\rho^2{v_kx_2\over u_kx_1}\,\delta_{k-1}\le{v_kx_2\over u_kx_1}\,
\rho^{a_1+\dots+a_{k-1}}/\rho^2\le K'\rho^{a_1+\dots+a_k},
\end{equation}
because $\displaystyle {v_k/ u_k}\le{\alpha^{a_k}/\alpha}$ and
$\pmatrix{x_1\cr x_2}=X_j$ is either $\pmatrix{1\cr \alpha}$ or
$\pmatrix{1\cr 0}$. The claim now follows from (\ref{comp3}),
(\ref{comp4}) and the fact that $\alpha^{(a_1+\cdots+a_k)/2}
\le\rho^n$.

\quad$\bullet$ {\it The case $\alpha>1$ with $\bfm=2$. -- } If
$\widehat\kappa+k$ is even, then
\begin{equation}\label{comp5}
|A_k|\le\rho^2\delta_k\le \rho^{a_1+\dots+a_k}/\rho^2.
\end{equation}
If $\widehat\kappa+k$ is odd, then
\begin{equation}\label{comp6}
|A_k|\le\rho^2{v_kx_2\over u_kx_1}\,\delta_{k-1}\le{v_k\over u_k}\cdot{x_2\over
x_1}\,
\rho^{a_1+\dots+a_k}/\rho^2,
\end{equation}
as $ {v_k/ u_k}\le{\alpha^{-a_k}/\alpha^{-1}}$ and
$\pmatrix{x_1\cr x_2}=\Delta X_j$ is necessarily
$\pmatrix{\alpha\cr 1}$ (because $\bfm=2$). The claim follows from
(\ref{comp5}), (\ref{comp6}) and the fact that
$\rho^{a_1+\cdots+a_k} \le\rho^n$.

\quad$\bullet$ {\it The case $\alpha=1$. -- } By (\ref{differences}),
\begin{eqnarray*}
\Big\vert\Phi(\xi)-\phi^w\Big\vert&\le&
\left\vert\log\left({p_{k_n}\over q_{k_n}}\right)-\log\left({p_{k_n-1}\over
q_{k_n-1}}\right)\right\vert+\\
&&\left\vert\left\{\log\left({p_{k_n}\over
q_{k_n}}\right)-\log\left({p_{k_n-1}\over q_{k_n-1}}\right)\right\}-
\left\{\log\left({p'_{k_n}\over q'_{k_n}}\right)-\log\left({p_{k_n-1}\over
q_{k_n-1}}\right)\right\}\right\vert,
\end{eqnarray*}
where ${p'_{k_n}/ q'_{k_n}}$ denotes the continued fraction
$[1;a_1,\dots,a_{k_n-1},a'_{k_n}]$. Since these continued
fractions are regular, it follows from the well known relations
(see \cite{Kin}) that
$$
\left\vert\Phi(\xi)-\phi^w\right\vert\le{2\over q_{k_n}q_{k_n-1}}+{1\over
q'_{k_n}q_{k_n-1}}.
$$
By induction, $q_{k_n}\ge a_1+\dots+a_{k_n}$ and $q'_{k_n}\ge
a_1+\cdots+a_{k_n-1}+a'_{k_n}=n$, whence
$\displaystyle\left\vert\Phi(\xi)-\phi^w\right\vert\le{3/
n}$.\hfill\break \null\hfill\qed
\medskip



\end{section}

\begin{section}{The multifractal analysis of the measure $\mu$}

{\bf 4.1.-- General case -- } The measure $\mu_*$ is $\gotF$-weak
Gibbs in all the cases described in Theorem~\ref{thA}, whence, in
view of Theorem~\ref{Fengolivier}, the following claim holds:

\begin{theorem}\label{thmultana0} {\sl Let $\mu$ be the $(\pp\;\qq)$-distributed $(\beta,2)$-Bernoulli
convolution, where $\beta$ is the multinacci number of degree
$\bfm\ge2$, with an extra condition $\pp\ge\qq$ if $\bfm\ge3$. The
multifractal domain $\hbox{\sc Dom}(\mu_*)$ is a compact interval
$[\underline{\alpha},\overline{\alpha}]$, where
$$
-\infty<\underline{\alpha}:=\lim_{q\to+\infty}{\tau_{\mu_*}(q)\over
q}\le \lim_{q\to-\infty}{\tau_{\mu_*}(q)\over
q}=:\overline{\alpha},
$$
and for any $\underline{\alpha}\leq \alpha\le\overline{\alpha}$,
$$
\dim_H
E_{\mu_*}(\alpha)=\sup_{q\in{\bf R}}\{\alpha q-\tau_{\mu_*}(q)\}.
$$}
\end{theorem}

We are now in position to state a multifractal formalism satisfied
by the Bernoulli convolution $\mu$ itself.

\begin{theorem}\label{thmultana}
{\sl Assume $\pp\ge\qq$ whenever $\bfm\ge3$; then
$\underline{\alpha}<\alpha<\overline\alpha$ whenever
$\dim_HE_\mu(\alpha)>0$, and for any $\underline{\alpha}\le
\alpha\le \overline{\alpha}$,
$$
\dim_HE_{\mu}(\alpha)=\sup_{q\in{\bf R}}\{\alpha q-\tau_{\mu_*}(q)\}.
$$
}
\end{theorem}

\begin{remark}{ (1) : We would like to emphasize that there exist cases
when $\tau_{\mu}(q)$is different from $\tau_{\mu_*}(q)$. This has
to do with the fact that the multifractal formalism in
Theorem~\ref{thmultana} is not complete in the sense that
$\hbox{\sc Dom}(\mu)$ may differ from the compact interval
$[\underline{\alpha};\overline{\alpha}]=\hbox{\sc Dom}(\mu_*)$
(see Lemma~\ref{nonconnex}).

\qquad (2) : Let $\mu'$ denote the  $(\qq\;\pp)$-distributed
$(\beta,2)$-Bernoulli convolution. Here one has
$\mu=\mu'\circ\SS$, where $\SS$ is the symmetry:
$\SS(x)=\alpha_\mu-x$. This clearly implies that, for any
$\alpha\in{\bf R}$
$$
\dim_H E_\mu(\alpha)=\dim_H E_{\mu'}(\alpha).
$$
Therefore, when $\pp<\qq$ and $\bfm\ge3$, the multifractal
formalism of $\mu$ is deduced by an application of
Theorem~\ref{thmultana} to the measure $\mu'$. }
\end{remark}

\qquad It remains to prove Theorem~\ref{thmultana}; actually, it
is a consequence of Theorem~\ref{thmultana0} and of the next
proposition. Loosely speaking, the latter asserts that $\mu$ has a
local Gibbs structure  whenever $\mu_*$ has a global one (recall
that the support of the measure $\mu$ is the interval
$[0,\alpha_\mu]$ with $\alpha_\mu=1/(\beta-1)$).

\begin{proposition}\label{comparaison}
{\sl For any $x$ in the support of $\mu$, there exists a constant
$K_x$ such that \hfil\break \null\qquad (i) : for $r$ small
enough, $\displaystyle{1\over K_x}\le {\mu(B_r(x))\over
\mu_*(B_r(x))}\le K_x$, if $x\in]0,1[;$ \hfil\break \null\qquad
(ii) : for $r$ small enough, $\displaystyle{1\over K_x}\le
{\mu(B_r(x))\over \mu_*(B_r(x-\ii_2))}\le K_x$, if
$\;x\in]1;\alpha_\mu[.$ }
\end{proposition}


{\bf  Proof. } (i) : Using the fact that $\mu$ and $\mu_0$
coincide on $[0,1[$, we first compare their values on the
$\gotF$-cylinders different from $\[\0^k\]$ for any $k\ge1$. To do
this, let $w=\0^k\eta w'$, for $k\ge 0$, $\eta\not=\0$ and $w'$ a
(possibly empty) word in ${\cal J}^*$; since
$\pmatrix{1&0}P_0=\pp^{\bfm}\pmatrix{1&0}$ and
$\pmatrix{1&1}P_0\le \pmatrix{1&1}$, we have
$$
{\mu\[w\]\over \mu_*\[w\]}
\ge{\pmatrix{1&0}P_\0^kP_\eta
P_{w'}R\over
\pmatrix{1&1}P_\0^kP_\eta
P_{w'}R} \ge \pp^{\bfm k}\cdot
{\pmatrix{1&0}P_\eta P_{w'}R\over
\pmatrix{1&1}P_\eta P_{w'}R}\;.
$$
Moreover, $\pmatrix{1&0}P_\eta \ge
\pp^\bfm\qq^\bfm\pmatrix{1&1}$ and
$\pmatrix{1&1}P_\eta
\le\pmatrix{1&1}$, whence
$$
{\mu\[w\]\over \mu_*\[w\]}
\ge \pp^{\bfm (k+1)}\cdot\qq^\bfm\cdot
{\pmatrix{1&1}P_{w'}V_\rho\over
\pmatrix{1&1}P_{w'}V_\rho}=\pp^{\bfm (k+1)}\cdot\qq^\bfm.
$$
Since the upper bound $\mu\[w\]/ \mu_*\[w\]\leq 2$ is always
valid, we obtain
\begin{equation}\label{Stilgar}
C^{k+1}\le{\mu\[w\]\over \mu_*\[w\]}\le 2,
\end{equation}
where $C=(\pp\qq)^\bfm$. Let $0<x<1$ and let $k_x\ge1$ be  an
integer such that $x\notin\[\0^{k_x}\]$; any ball $B_r(x)\subset\,
]0,1[$ which does not intersect $\[\0^{k_x+1}\]$ can be tiled by
countably many cylinders of the form $\[\0^k\eta w'\]$, where
$k\leq k_x$, $w\not=0$ and $w'\in{\cal J}^*$. Now it follows from
(\ref{Stilgar}) that
\begin{equation}\label{Stilgar'}
C^{k_x+1}\leq
{\mu(B_r(x))\over
\mu_*(B_r(x))}\leq 2\;.
\end{equation}

\qquad (ii) : Now, let $x\in\left]1,\alpha_\mu\right[$. The
Bernoulli convolutions $\mu$ and $\mu'$ associated with the
probability vectors $(\pp,\qq)$ and $(\qq,\pp)$ respectively,
satisfy the relation $\mu=\mu'\circ\SS$, where
$\SS(t)=\alpha_\mu-t$ for any $t\in[0,\alpha_\mu]$. Since
$x'=\SS(x)$ belongs to $]0,1[$, there exists $k_x\ge1$ such that
$x'\notin\[\0^{k_x}\]$; we apply (\ref{Stilgar'}) to the measures
$\mu'$  and $\mu'_*$ at $x'$ and obtain
$$
C^{k_x+1}\leq
{\mu'(B_r(x'))\over
\mu'_*(B_r(x'))}\leq 2\;.
$$
This yields (ii), as $\mu'(B_r(x'))=\mu(B_r(x))$, and
$$
\displaystyle\mu'_*(B_r(x'))=
{\mu'(B_r(x'))+\mu'\big(B_r(x'+\ii_2)\big)\over\mu'([0\,;1])+\mu'([0\,;1]+\ii_2)}=
\lambda\mu_*(B_r(x-\ii_2)),
$$
where $\lambda>0$ is a constant.\hfill\break
\null\hfill\qed
\medskip

\qquad {\bf 4.2.-- Case of the Erd\H os measure --} As we have
seen, ``local" Gibbs properties of the Bernoulli  convolution in a
multinacci base are sufficient to establish the multifractal
formalism of the level sets with positive Hausdorff dimension.
However, a natural question would be to determine whether or not
there exists a reasonable net with respect to which the measure in
question has a ``global" Gibbs structure.

\qquad In the rest of the paper we concentrate on the case of the
{\it Erd\H os measure}, {\it i.e.} $\bfm=2$ or, equivalently,
$\beta=(1+\sqrt{5})/2$. Following our notation introduced in
Section~3.1, the $\beta$-shift is associated with the two affine
contractions:
$$
\RR_0(x)=x/\beta\quad\hbox{and}\quad\RR_1(x)=x/\beta+1/\beta.
$$
The corresponding s.a.c. ${\cal R}$ consists of the three
contractions:
$$
\widehat{\RR}_\0(x)\!=\!\RR_{00}(x)\!=\!x/\beta^2,\;
\widehat{\RR}_\1(x)\!=\!\RR_{010}(x)\!=\!x/\beta^3\!+\!1/\beta^2,\;
\widehat{\RR}_\2(x)\!=\!\RR_{10}(x)\!=\!x/\beta^2\!+\!1/\beta.
$$
It is adapted to the interval $[0,1]$ (see {\it Figure~4}). Recall
that the measure $\mu$ is supported by the interval $[0,\beta]$;
to study the global Gibbs properties of $\mu$, it is thus more
convenient to make an affine scale change from $[0,1]$ to
$[0,\beta]$. Thus, instead of ${\cal R}$ we consider the s.a.c.
${\cal S}=\{\SS_\0,\SS_\1,\SS_\2\}$ with
$$
\SS_\0(x)=x/\beta^2,\quad
\SS_\1(x)=x/\beta^3+1/\beta,\quad\SS_\2(x)=x/\beta^2+1.
$$
Clearly, ${\cal S}$ is adapted to the interval $[0;\beta]$ and we
denote the associated $3$-fold net by $\widetilde{\gotF}$. For any
word $w\in\{\0,\1,\2\}^*$, we put $\ltriple
w\rtriple=\SS_w[0,\beta[$; obviously, the set of $\ltriple
w\rtriple$ determines the basic intervals of $\widetilde{\gotF}$.
Moreover, any Borel measure $\nu$ on the real line is associated
with the measure $\widetilde{\nu}$ defined on an arbitrary
interval~$J$ as follows: $\widetilde{\nu}(J)=\nu(J/\beta)$. This
scale change clearly implies that the measures
$\widetilde{\mu}_0$, $\widetilde{\mu}_2$ and $\widetilde{\mu}_*$
satisfy, for any word $w\in\{\0,\1,\2\}^*$, the following matrix
identities:
$$
\pmatrix{\widetilde{\mu}_0\ltriple
w\rtriple\cr\widetilde{\mu}_2\ltriple
w\rtriple}=P_w\pmatrix{\mu_0([0,1])\cr\mu_2([0,1])}
\quad\hbox{and}\quad \widetilde{\mu}_*\ltriple
w\rtriple=\pmatrix{1&1}P_w\pmatrix{\pp\cr\qq},
$$
where
$$
P_\0=\pp^2\pmatrix{1&0\cr\qq/\pp&\qq/\pp},\;\;
P_\1=\pp^2\qq\pmatrix{1&1\cr\qq/\pp&\qq/\pp},\;\;
P_\2=\qq^2\pmatrix{\pp/\qq&\pp/\qq\cr0&1}.
$$
We are going to use two key properties of this model. Firstly, the
probability measure $\mu$ satisfies the following well known
self-similar equation:
$$
\mu=\pp\mu\circ\RR_0^{-1}+\qq\mu\circ\RR_1^{-1}.
$$
Secondly,
$$
\SS_\0=\RR_{00},\quad
\SS_\1=\RR_{100}=\RR_{011},\quad\SS_\2=\RR_{11},
$$
where the identity $\RR_{100}=\RR_{011}$ plays a crucial role. Let
$J$ be a subinterval of $[0,\beta]$; then, one has successively
\begin{eqnarray*}
\bullet\;\;\mu\Big(\SS_\0(J)\Big) & = &
\pp\mu\Big(\RR_0^{-1}\RR_{00}(J)\Big)+\qq\mu\Big(\RR_1^{-1}\RR_{00}(J)\Big)\\
                                  & = &
\pp\mu\Big(\RR_{0}(J)\Big)=\pp\widetilde{\mu}_0(J)\;;
\\\\
\bullet\;\;\mu\Big(\SS_\1(J)\Big) & = &
\pp\mu\Big(\RR_0^{-1}\RR_{011}(J)\Big)+\qq\mu\Big(\RR_1^{-1}\RR_{100}(J)\Big)\\
                                  & = &
\pp\mu\Big(\RR_{11}(J)\Big)+\qq\mu\Big(\RR_{00}(J)\Big)\\
                                  & = &
\pp\qq\mu\Big(\RR_{1}(J)\Big)+\qq\pp\mu\Big(\RR_{0}(J)\Big)=
\pp\qq\Big(\widetilde{\mu}_0(J)+\widetilde\mu_2(J)\Big)\;;
\\\\
\bullet\;\;\mu\Big(\SS_\2(J)\Big) & = &
\pp\mu\Big(\RR_0^{-1}\RR_{11}(J)\Big)+\qq\mu\Big(\RR_1^{-1}\RR_{11}(J)\Big)\\
                                  & = &
\qq\mu\Big(\RR_{1}(J)\Big)=\qq\widetilde{\mu}_2(J).
\end{eqnarray*}
We conclude by the fact that, for any $\eta\in\{\0,\1,\2\}$ and any
$w\in\{\0,\1,\2\}^*$, one has
\begin{equation}\label{defmu}
\mu\ltriple\eta w\rtriple=V_\eta P_{w}
\pmatrix{\displaystyle{\pp/(1-\pp\qq)}\cr \displaystyle{\qq/(1-\pp\qq)}},
\quad\hbox{where}\quad
\cases{
V_\0=\pmatrix{\pp&0}\;;\cr
V_\1=\pmatrix{\pp\qq&\pp\qq}\;;\cr
V_\2=\pmatrix{0&\qq}.}
\end{equation}

Consider two subcases.

\qquad {\bf 4.2.1.-- The uniform case --} We first consider the
uniform Erd\H os measure, {\em i.e.}, $\pp=\qq=1/2$. Then for
$\eta\in\{\0,\1,\2\}$ and $w\in\{\0,\1,\2\}^*$:
\begin{equation}\label{****}
\mu\ltriple\eta w\rtriple=V_\eta' P_{w}
\pmatrix{1\cr 1},
\;\hbox{where}\;
\cases{
V_\0'=\pmatrix{1/3&0}\cr
V_\1'=\pmatrix{1/6&1/6}\cr
V_\2'=\pmatrix{0&1/3}}
\;\hbox{and}\;
\cases{P_\0=\displaystyle {1\over 4}\pmatrix{1&0\cr1&1}\cr
P_\1=\displaystyle{1\over 4}\pmatrix{1/2&1/2\cr1/2&1/2}\cr
P_\2=\displaystyle{1\over 4}\pmatrix{1&1\cr0&1}.\cr}
\end{equation}

\qquad We are going to show that $\mu$ is $\widetilde{\gotF}$-weak
Gibbs. Consider the probability measure $\widetilde{\mu}_*$ with
the support equal to the interval $[0,\beta]$ defined as follows:
for any word $w\in\big\{\0,\1,\2\big\}^*$,
\begin{equation}\label{BC6}
\widetilde{\mu}_*\ltriple w\rtriple= {1\over 2} \pmatrix{1&1}
P_w\pmatrix{1\cr1}.
\end{equation}
By Theorem~\ref{Holder}, $\widetilde{\mu}_*$ is a
$\widetilde{\gotF}$-weak Gibbs measure with the potential
$\Phi:\{\0,\1,\2\}^{\bf N}\to{\bf R}$; the formula for $\Phi$ can
in fact be given by means of regular continued fractions: let
$a_0,\dots,a_{2n}$ be $2n$ integers ($n\ge0$) with
$a_0,\dots,a_{2n-1}>0$ and $a_{2n}\ge0$, when $n\ge1$. Then for
any $\xi\in\{\0,\1,\2\}^{\bf N}$,
$$
\Phi(\0^{a_0}\2^{a_1}\cdots
\2^{a_{2n-1}}\0^{a_{2n}}\1\xi)=\Phi(\2^{a_0}\0^{a_1}\cdots
\0^{a_{2n-1}}\2^{a_{2n}}\1\xi)=\log\Big({f(a_0,\dots,a_{2n})/
4}\Big)
$$
with $f(0)=1$ and
$$
f(a_0,\dots,a_{2n})=1+
{{1}\over{\displaystyle a_0+
{{1\raise10pt\hbox{}}\over{\displaystyle\;\;\,\raise-8pt\hbox{$\ddots$}\;\;\,
\raise-20pt\hbox{$\displaystyle +
{1\over a_{2n}+1}$}
\!\!\!\!\!\!\!\!\!\!\!\!}}}
\!\!\!\!\!\!\!\!\!}\qquad\qquad.
$$

\begin{theorem}\label{erdosunif}
{\sl The uniform Erd\H os measure $\mu$ is a
$\widetilde{\gotF}$-weak Gibbs measure of the potential~$\Phi$.}
\end{theorem}

This theorem is a consequence of the fact that $\widetilde{\mu}_*$
is itself a $\widetilde{\gotF}$-weak Gibbs measure of $\Phi$ and
of the following proposition:

\begin{proposition}\label{Erdpr1}{\sl For any $\omega\in
\{\0,\1,\2\}^{\bf N}$ and any integer $n\ge1$,
$$
{8\over
3(n+2)}\leq{\mu\ltriple\omega_0\cdots\omega_{n-1}\rtriple\over
\widetilde{\mu}_*\ltriple\omega_0\cdots\omega_{n-1}\rtriple}\leq
{8\over 3}.
$$}
\end{proposition}

{\bf Proof.} Note first that $\mu\ltriple
w\rtriple/\widetilde{\mu}_*\ltriple w\rtriple=4/3$, whenever
$w=\1w'$; moreover,  if $w=\0^n$ or $\2^n$ with $n\ge1$, then
$\mu\ltriple w\rtriple/\widetilde{\mu}_*\ltriple
w\rtriple=8/(3(n+2))$. Now, given $n\ge2$, we assume that
$w=\eta^a\nu w'\in\{\0,\1,\2\}^n$, with $a<n$. Without loss of
generality we assume $\eta=\0$ and $\nu=\2$ (the other cases with
$\eta\in\{\0,\2\}$ and $\eta\ne\nu\in\{\0,\1,\2\}$ are similar).
It is straightforward that \hbox{$\mu\ltriple
w\rtriple/\widetilde{\mu}_*\ltriple w\rtriple\leq 8/3$}, so it
remains to establish the lower bound of $\mu\ltriple
w\rtriple/\widetilde{\mu}_*\ltriple w\rtriple$. As
$$
P_{w'} \pmatrix{1\cr 1}:=\pmatrix{x\cr y},
$$
we have
$$
{\mu\ltriple w\rtriple\over \widetilde{\mu}_*\ltriple w\rtriple}
={8\over 3\cdot4^{a+1}}\cdot{{\pmatrix{1&1} \pmatrix{x\cr y}}\over
      {\pmatrix{1&1} P_0^{a}\,P_2\pmatrix{x\cr y}}}=
{8\over 3}\cdot{(x+y)\over (1+a)(x+y)+y}\geq{8\over
    3(a+2)}.
$$
Since $a<n$, we conclude that $\mu\ltriple
w\rtriple/\widetilde{\mu}_*\ltriple w\rtriple\geq{8\over
3(n+2)}$.\hfill\break \null\hfill \qed
\medskip

\qquad We would like to stress that $\mu$ is {\em not} a
$\widetilde{\gotF}$-Gibbs measure. Actually, if
\hbox{$\overline0:=(\omega_i=0)_{i=0}^\infty$} then, for any
potential $\psi:\Sigma\to{\bf R}$, we have
$$
\exp(S_n\psi(\2\overline\0))={\exp(\psi(\2\overline\0))\over
\exp(\psi(\overline\0))}\Big\{\exp(\psi(\overline\0))\Big\}^n.
$$
A direct computation, in view of (\ref{****}), yields
$\mu\ltriple\2\0^{n-1}\rtriple=n/(3\cdot 4^{n-1})$, for any $n>0$;
if $\mu$ were a Gibbs measure of $\psi$, then there would exist a
constant $K>1$ such that for any $n\ge1$,
$$
{1 \over K}\leq{1\over n}\cdot
\Big\{4\cdot\exp\big(\psi(\overline0)\big)\Big\}^n\le
K,
$$
which is impossible. Hence $\mu$ is not $\widetilde{\gotF}$-Gibbs.

\qquad Furthermore, in this specific case one can prove a much
more refined result, namely that $\mu$ is not Gibbs in a very
strong sense---see Theorem \ref{local2} below.

\qquad Assume that $\gotE$ is an arbitrary ${\bf s}$-fold net
which is quasi-Bernoulli; if $\mu$ is $\gotE$-Gibbs then, in
particular, it is $\gotE$-quasi-Bernoulli and by
Theorem~\ref{herteau}~(i), its scale spectrum $\tau_{\mu}$ is
differentiable on the whole real line. However, this is not the
case, since it is known \cite{Fen, FO} that there exists $q_c<-2$
such that $\tau_{\mu}$ is not differentiable at $q_c$.

Thus, we have proved

\begin{theorem}\label{local2}
{\sl  There is no quasi-Bernoulli net with respect to which the uniform
Erd\H os measure  is quasi-Bernoulli.}
\end{theorem}

\qquad {\bf 4.2.2.-- The nonuniform case --} From now on we assume
$\pp\ne\qq$. As the measure $\widetilde{\mu}_*$ is
$\widetilde{\gotF}$-Gibbs, its multifractal domain is a compact
interval. In fact,
$$
\hbox{\sc Dom}({\mu}_*)=\hbox{\sc
Dom}(\widetilde{\mu}_*):=[\underline{\alpha},\overline\alpha],
$$
where
$$
\underline{\alpha}=
\inf_{x\in[0,\beta]}\left\{\lim_{r\to0}{\log\widetilde{\mu}_*(B_r(x
))\over \log(r)}\right\}= \inf_{\omega\in\{\0,\1,\2\}^{\bf
N}}\left\{\lim_{n\to\infty}{\log\widetilde{\mu}_*\ltriple\omega_0\cdots
\omega_{n-1}\rtriple\over
\log\vert\ltriple\omega_0\cdots\omega_{n-1}\rtriple\vert}\right\}
$$
and
$$
\overline{\alpha}=
\sup_{x\in[0,\beta]}\left\{\lim_{r\to0}{\log\widetilde{\mu}_*(B_r(x)
)\over \log(r)}\right\}= \sup_{\omega\in\{\0,\1,\2\}^{\bf
N}}\left\{\lim_{n\to\infty}{\log\widetilde{\mu}_*\ltriple\omega_0\cdots
\omega_{n-1}\rtriple\over
\log\vert\ltriple\omega_0\cdots\omega_{n-1}\rtriple\vert}\right\}.
$$

\begin{theorem}\label{local3}
{\sl  Assume that $\pp\ne\qq$; then there is no quasi-Bernoulli
net with respect to which the corresponding (nonuniform) Erd\H os
measure is quasi-Bernoulli or weak Gibbs.}
\end{theorem}

If $\mu$ were a weak Gibbs measure w.r.t. a quasi-Bernoulli net of
the interval $[0,\beta]$, it would be necessary that its
multifractal domain were a compact interval. Moreover, by
Proposition~\ref{comparaison}, one would have
$$
\hbox{\sc Dom}(\mu)=\hbox{\sc Dom}({\mu}_*)=\hbox{\sc
Dom}({\widetilde\mu}_*)=[\underline\alpha,\overline\alpha].
$$
Thus, Theorem~\ref{local3} would follow from

\begin{lemma}\label{nonconnex}
{\sl  If $\pp\ne\qq$ then
$$
\hbox{\sc Dom}(\mu)\not=\hbox{\sc Dom}({\mu}_*)=\hbox{\sc
Dom}({\widetilde\mu}_*)=[\underline\alpha,\overline\alpha].
$$}
\end{lemma}

{\bf Proof.} Assume without loss of generality that
$\pp<\qq$---the case $\pp>\qq$ is handled by considering the
symmetry between the $(\pp,\qq)$ and the $(\qq\;\pp)$-distributed
Erd\H os measure. On one hand, given any
$\omega\in\{\0,\1,\2\}^{\bf N}$, as is easy too see, there exists
a constant $K>0$ such that
$$
\widetilde{\mu}_*\ltriple\omega_0\cdots\omega_{n-1}\rtriple\ge
K(\pp\qq)^{n}
$$
so that
$$
{\log\widetilde{\mu}_*\ltriple\omega_0\cdots\omega_{n-1}\rtriple\over
\log\vert\ltriple\omega_0\cdots\omega_{n-1}\rtriple\vert}\le {\log
K+n\log(\pp\qq)\over 2n\log(1/\beta)},
$$
whence
\begin{equation}\label{bidule2}
\overline\alpha\leq {\log(\pp\qq)\over 2\log(1/\beta)}.
\end{equation}
On the other hand, a direct computation yields that for any integer $n\ge0$,
$$
\mu\ltriple\0^n\rtriple=
\pp^{2n-2}\pmatrix{\pp&0}\pmatrix{1&0\cr\qq/\pp&\qq/\pp}^{n-1}\pmatrix{\pp
/(1-\pp\qq)\cr\qq/(1-\pp\qq)}=K'\pp^{2n}
$$
so that
\begin{equation}\label{bidule3}
\lim_{r\to0}{\log\mu\ltriple\0^{n}\rtriple\over
\log\vert\ltriple\0^{n}\rtriple\vert}={\log\pp\over \log(1/\beta)}
\end{equation}

Given any $0<r<\beta$, let $n_r$ be the integer satisfying
$1/\beta^{2(n_r+1)}\le r<1/\beta^{2n_r}$, whence
$$
{n_r\log(1/\beta^2)\over \log
r}\cdot{\log\mu\ltriple\0^{n_r}\rtriple\over
\log\vert\ltriple\0^{n_r}\rtriple\vert}\le{\log\mu(B_r(0))\over\log
r}\le {(n_r+1)\log(1/\beta^2)\over\log
r}\cdot{\log\mu\ltriple\0^{(n_r+1)}\rtriple\over
\log\vert\ltriple\0^{(n_r+1}\rtriple\vert}.
$$
Since $\pp<\qq$, we have, in view of (\ref{bidule2}) and
(\ref{bidule3}),
$$
\lim_{r\to0}{\log\mu(B_r(0))\over \log r}={\log\pp\over
\log(1/\beta)}>{\log(\pp\qq)\over
2\log(1/\beta)}\ge\overline\alpha.
$$
\hfill\break
\null\hfill\qed

\begin{remark}{One can prove that the multifractal domains of $\mu$ and $\mu_*$ differ because of the
local dimension of $\mu$ at the point $x=0$ if $\pp<\qq$ and
$x=\alpha_\mu$ if $\pp>\qq$. Therefore, the multifractal domain of
the nonuniform Erd\H os measure is the disconnected union of an
interval and a singleton. More precisely,
$$
\hbox{\sc
Dom}(\mu)=[\underline{\alpha},\overline\alpha]\cup\{\alpha^*\}\quad\hbox{with}\quad
\alpha^*=\max\left\{{\log\pp\over
\log(1/\beta)},{\log\qq\over
\log(1/\beta)}\right\}
$$

\qquad A similar pattern holds in the case of the $3$-fold
convolution of the Cantor measure. Namely, let $\gamma$ denote the
Cantor measure, {\em i.e.}, the self-similar probability measure
associated with the two affine contractions
$$
\SS_{\bf 0}:x\mapsto x/3\quad\hbox{and}\quad\SS_{\bf 1}:x\mapsto x/3+2/3.
$$
It can be easily checked that the $3$-fold convolution measure
$\gamma':=\gamma*\gamma*\gamma$ is none other than the
$\pp$-distributed $(2,3)$-Bernoulli convolution with
$$
\pp=\pmatrix{\displaystyle{1\over8}& \displaystyle{3\over
8}&\displaystyle{3\over 8}&\displaystyle{1\over 8}\cr}.
$$
A detailed multifractal analysis of $\gamma'$ presented in
\cite{HL} shows that the multifractal domain $\hbox{\sc
Dom}(\gamma')$ is no longer a compact interval. Actually,
$\hbox{\sc Dom}(\gamma')$ is proved to be the union of a compact
interval  and a singleton. This proves that $\gamma'$ is neither
quasi-Bernoulli nor weak Gibbs w.r.t. any given reasonable ${\bf
s}$-fold net of the unit interval. }
\end{remark}

\medskip

\vbox{\qquad{\bf Acknowledgment -- }{\sl A part of this paper was
written during the postdoctoral visit of the first author at the
Department of Mathematics of the Chinese University of Hong Kong.
He is grateful to Jean-Pierre Kahane and Ka-Sing~Lau for their
constant support and valuable discussions.} }

\end{section}
\vfill
\eject

{\baselineskip=12pt



\baselineskip=14pt

Eric OLIVIER $<$Eric.Olivier@up.univ-mrs.fr$>$\par
SCAM, Universit\'e de Provence,\par
3, place Victor Hugo,\par
13331 Marseille Cedex 3, France.
\medskip

Nikita SIDOROV $<$sidorov@manchester.ac.uk$>$\par School of
Mathematics, The University of Manchester, P.O. Box 88, \par
Manchester M60 1QD, United Kingdom.
\medskip

Alain THOMAS $<$thomas@gyptis.univ-mrs.fr$>$\par
Centre de
Math\'ematiques et d'Informatique, LATP \'Equipe DSA,\par
9, rue F. Joliot-Curie,\par
13453 Marseille Cedex 13, France.

}

\end{document}